\documentclass[10pt,a4paper,leqno]{amsart}

\usepackage{amsmath,amssymb,amsthm}

\def\dis
{\displaystyle}
\def\eps
{\varepsilon}

\def\R{{\mathbb R}}
\def\N{{\mathbb N}}
\def\Z{{\mathbb Z}}

\def\F{{\mathcal F}}
\def\S{{\mathcal S}}
\def\virgp{\raise 2pt\hbox{,}}

\def\Eq#1#2{\mathop{\sim}\limits_{#1\rightarrow#2}}
\def\Tend#1#2{\mathop{\longrightarrow}\limits_{#1\rightarrow#2}}
\def\Norm#1#2{\left\| #1 \right\|_{#2}}

\def\d{{\partial}}
\def\e{\varepsilon}
\def\si{{\sigma}}
\def\U{{\tt U}}
\def\V{{\tt V}}
\def\W{{\tt W}}

\def\v{{\tt v}}
\def\w{{\tt w}}
\def\H{{\tt H}}
\def\G{{\tt G}}
\def\O{\mathcal O}
\def\<{\langle}
\def\>{\rangle}

\def\1{{\rm 1\mskip-4.5mu l} }

\def\({\left(}
\def\){\right)}
\def\supp{\operatorname{supp}}

\theoremstyle{plain}
\newtheorem{theo}{Theorem}[section]
\newtheorem{lem}[theo]{Lemma}
\newtheorem{cor}[theo]{Corollary}
\newtheorem{prop}[theo]{Proposition}

\theoremstyle{definition}
\newtheorem{defin}[theo]{Definition}
\newtheorem*{nota}{Notation}

\theoremstyle{remark}
\newtheorem{rema}[theo]{Remark}
\newtheorem*{rema*}{Remark}

\numberwithin{equation}{section}

\begin{document}

\title[Quadratic oscillations in NLS II. The $L^2$-critical
case]{On the role of 
quadratic oscillations in nonlinear Schr{\"o}dinger equations II. The
$L^2$-critical case} 
\author[R. Carles]{R{\'e}mi Carles}
\address{IRMAR, Universit\'e de Rennes 1\\ Campus de
Beaulieu\\ 35~042 Rennes cedex\\ France}
\email{remi.carles@math.univ-rennes1.fr}
\author[S. Keraani]{Sahbi Keraani}
\email{sahbi.keraani@univ-rennes1.fr}
\thanks{
This work was partially supported by the ACI grant
  ``{\'E}quation des ondes : oscillations, dispersion et
contr{\^o}le'', and by the European network HYKE, funded by the EC as
contract HPRN-CT-2002-00282.}
\begin{abstract}
We consider a nonlinear semi--classical Schr{\"o}dinger equation
for which quadratic oscillations lead to focusing at
one point, described by a nonlinear scattering operator. The relevance
of the nonlinearity was discussed by R.~Carles, C.~Fermanian--Kammerer
and I.~Gallagher for
$L^2$-supercritical power-like nonlinearities and more general initial
data. The present results
concern the $L^2$-critical case, in space dimensions $1$ and $2$; we
describe the set of 
non-linearizable data, which is larger, due to the scaling. As
an application, we precise a result by F.~Merle and L.~Vega concerning finite
time blow up  for the critical Schr{\"o}dinger equation. The proof
relies on linear and nonlinear profile decompositions.
\end{abstract}
\maketitle

\section{Introduction}
\label{sec:intro}

Consider the initial value problem
\begin{equation}\label{eq:cfg}
i\e \d_t u^\e +\frac{1}{2}\e^2\Delta u^\e = \e^{n\si}|u^\e|^{2\si}u^\e\
;\quad u^\e_{\mid t=0} =u_0^\e\, ,
\end{equation}
where $x\in \R^n$ and $\e\in ]0,1]$. Our aim is to understand the
relevance of the nonlinearity in the limit $\e \to 0$, according to
the properties of the initial data $u_0^\e$. In
\cite{CFG}, the case $\si >2/n$, with $\si <2/(n-2)$ if 
$n\geq 3$ and $u_0^\e$, $\e
\nabla_x u_0^\e$ bounded in $L^2(\R^n)$ uniformly for $\e \in
]0,1]$, was studied. Note that under these assumptions, global
existence in $H^1(\R^n)$ for fixed $\e >0$ is well known (see
e.g. \cite{CazCourant}). It was proven that the nonlinearity has a leading
order influence in the limit $\e \to 0$ if and 
only if the initial data include a quadratic oscillation of the form
\begin{equation*}
f(x-x^\e) e^{-i\frac{|x-x^\e|^2}{2\e t^\e}}\, ,
\end{equation*}
for some $x^\e\in \R^n$ and $t^\e >0$, with $\limsup_{\e \to 0} t^\e
/\e \in ]0,+\infty[$ (see \cite[Theorem~1.2]{CFG}  for a precise
statement). Two things have to be said about this property. First, it
shows that the presence of quadratic oscillationsis necessary for the
nonlinearity to have a leading order influence; it was established in
\cite{Ca2} that it is sufficient.  Second,
only one scale is involved 
in such initial profiles, that is, $\e$. In the present paper, we
study the $L^2$-critical case, $\si =2/n$. We prove that quadratic
oscillations are not necessary to have a leading
order nonlinear behavior, if we assume that the initial data satisfy
the same assumptions as in \cite{CFG}; other scales than $\e$
have to be taken into account, because $\si =2/n$
corresponds to the critical scaling at the $L^2$ level. To see
this, consider a solution $\U$ to the nonlinear 
Schr\"odinger equation
\begin{equation}\label{eq:NLS}
i\d_t \U +\frac{1}{2}\Delta \U = \lambda |\U|^{4/n}\U\, ,
\end{equation}
with $\lambda =1$ and $\U_{\mid t=0}=\phi$.
If $\phi\in \Sigma$, where
\begin{equation}\label{eq:Sigma}
\Sigma := \left\{ \phi \in H^1(\R^n)\, ;\ |x|\phi \in
L^2(\R^n)\right\}\, ,
\end{equation}
then the solution $\U$ of \eqref{eq:NLS} is defined globally in
time, with $\U\in C(\R_t;\Sigma)$ (see e.g. \cite{CazCourant}). Let
$(t_0,x_0)\in \R \times \R^n$. It is straightforward
to see that 
\begin{equation}\label{eq:exsqrt}
u^\e (t,x) = \frac{1}{\e^{n/4}}\U \left(t-t_0, \frac{x-x_0}{\sqrt\e}\right)
\end{equation}
solves \eqref{eq:cfg} with $\si=2/n$, and that $u^\e(0,\cdot)$ and $\e
\nabla_xu^\e(0,\cdot)$ are bounded in $L^2(\R^n)$, uniformly in $\e
\in ]0,1]$. This particular solution is such that the nonlinearity in
\eqref{eq:cfg} is relevant at leading order, at any (finite) time, near
$x=x_0$. 
This is in contrast with the $L^2$-supercritical case $\si
>2/n$, where only profiles of the form 
\begin{equation}\label{eq:excfg}
u^\e (t,x) = \frac{1}{\e^{n/2}}\U \left(\frac{t-t^\e}{\e}\virgp 
\frac{x-x^\e}{\e}\right)
\end{equation}
were relevant. The solutions \eqref{eq:exsqrt} are deduced from the
solutions \eqref{eq:excfg} by scaling. If $\U$ solves \eqref{eq:NLS},
then so does $\widetilde\U$, given by 
\begin{equation*}
\widetilde\U(t,x) = \lambda^{n/2}
\U\left( \lambda^2 t, \lambda x\right) ,
\end{equation*}
for any real $\lambda$: the case
$\si=2/n$ is $L^2$-critical. Applying this transform to
solutions \eqref{eq:excfg} with $\lambda =\sqrt\e$
yields solutions \eqref{eq:exsqrt}, with $t_0 =t^\e/\e$ and $x_0 =x^\e
/\sqrt\e$. 

Before going further into details, we fix some notations and introduce a
definition. We consider initial value problems
\begin{equation}\label{eq:ck}
i\e \d_t u^\e +\frac{1}{2}\e^2\Delta u^\e = \lambda\e^2|u^\e|^{4/n}u^\e\
;\quad u^\e_{\mid t=0} =u_0^\e\, ,
\end{equation}
with $\lambda \in \{-1,+1\}$, that is, we consider the $L^2$-critical
case of \eqref{eq:cfg}, with possibly focusing nonlinearities
($\lambda =-1$). As in \cite{CFG}, we define the free
evolution $v^\e$ of $u_0^\e$, 
\begin{equation}\label{eq:v}
i\e \d_t v^\e +\frac{1}{2}\e^2\Delta v^\e = 0\
;\quad v^\e_{\mid t=0} =u_0^\e\, .
\end{equation}
We resume some notations used in \cite{CFG}.
\begin{nota}
i) For a family $(a^\eps)_{0<\eps \leq 1}$ of functions in
$H^1(\R^n)$, define
\begin{equation*}
\|a^\eps\|_{H^1_\eps}  :=  \|a^\eps\|_{L^2} + \|\eps \nabla
a^\eps\|_{L^2}\, .
\end{equation*} 
We will say that $a^\eps$ is
bounded (resp. goes to zero) in $H^1_\eps$ if 
$$
\limsup_{\eps \rightarrow 0 }\|a^\eps\|_{H^1_\eps}<\infty  
\textrm{ (resp. }=0\textrm{)}.$$
ii) If $(\alpha^\eps)_{0<\eps \leq 1}$ and $(\beta^\eps)_{0<\eps \leq 1}$
are two families of positive numbers, we write
\begin{equation}
\alpha^\eps \lesssim \beta^\eps
\end{equation}
if there exists $C$ independent of $\eps \in ]0,1]$ such that for
any $\eps \in ]0,1]$, 
$$\alpha^\eps \leq C \beta^\eps.$$
\end{nota}
From now on, $u^\e$ (resp. $v^\e$) stands for the solution to
\eqref{eq:ck} (resp. \eqref{eq:v}), with 
$\lambda=-1$ or $+1$ indifferently, unless precised specifically. 
\begin{defin}[Linearizability] Let $u_0^\e\in L^2(\R^n)$, 
bounded in $L^2(\R^n)$, and let $I^\e$ be an interval
of $\R$, possibly depending on $\e$. 
\begin{itemize}
\item[i)] The solution $u^\e$ is linearizable on $I^\e$ in $L^2$ if
\begin{equation*}
\limsup_{\e \to 0}\sup_{t\in I^\e}\left\|
u^\e(t)-v^\e(t)\right\|_{L^2(\R^n)} =0\, .
\end{equation*}
\item[ii)] If in addition $u_0^\e\in H^1(\R^n)$ and $u_0^\e$ is
bounded in $H^1_\e$, we say that $u^\e$ is linearizable on $I^\e$ in
$H^1_\e$ if 
\begin{equation*}
\limsup_{\e \to 0}\sup_{t\in I^\e}\left( \left\|
u^\e(t)-v^\e(t)\right\|_{L^2(\R^n)}+\left\|\e \nabla_x
u^\e(t)-\e\nabla_xv^\e(t)\right\|_{L^2(\R^n)}\right) =0\, .
\end{equation*}
\end{itemize}
\end{defin}
We prove the following result. Notice that we have to restrict to
the case of space dimensions $1$ and $2$ (see
Remark~\ref{rema:dimension} below). 
\begin{theo}\label{theo:cns} Assume $n=1$ or $2$. Let $u_0^\e$ 
bounded in $L^2(\R^n)$, $I^\e \ni 0$ a time interval. 
\begin{itemize}
\item $u^\e$ is linearizable on $I^\e$ in $L^2$ if and only if
\begin{equation}\label{eq:linear}
\limsup_{\e \to 0}\e \|v^\e\|_{L^{2+4/n}(I^\e\times \R^n)}^{2+4/n} = 0
\, .
\end{equation}
\item Assume in addition that $u_0^\e\in H^1$ and $u_0^\e$ is
bounded in $H^1_\e$. Then $u^\e$ is linearizable on $I^\e$ in $H^1_\e$ if
and only if \eqref{eq:linear} holds. 
\end{itemize}
\end{theo}
Notice that a similar result was proven in \cite{CFG}, in the
$L^2$-supercritical case, with a different linearizability condition,
\begin{equation}\label{eq:lincfg}
\limsup_{\e \to 0}\e^2\sup_{t\in I^\e} \|v^\e(t)\|_{L^{2+4/n}(
\R^n)}^{2+4/n} = 0 
\, .
\end{equation}
The fact that this condition is necessary for $u^\e$ to be
linearizable in $H^1_\e$ is easy to see, from the classical
conservations of mass and energy, which we write in the case $\si=2/n$
(in the general case, the powers $2+4/n$ are replaced by $2\si +2$):
\begin{equation}\label{eq:conservations}
\begin{aligned}
&\text{Mass: }
\frac{d}{dt}\|u^\e(t)\|_{L^2}=\frac{d}{dt}\|v^\e(t)\|_{L^2}=0 \, .\\
&\text{Linear energy: } \frac{d}{dt}\|\e \nabla_x v^\e(t)\|_{L^2}=0\, .\\
&\text{Nonlinear energy: } \frac{d}{dt}\left(\frac{1}{2}\|\e \nabla_x
u^\e(t)\|_{L^2}^2 
+\frac{\lambda \e^2}{2+4/n}\|u^\e(t)\|_{L^{2+4/n}}^{2+4/n}\right)=0\, .
\end{aligned}
\end{equation}
The proof that condition \eqref{eq:lincfg} implies linearizability in
$H^1_\e$ (which is a stronger property than linearizability in $L^2$)
involves Strichartz estimates, and seemed to rely in an unnatural way
on the assumption $\si >2/n$. Example \eqref{eq:exsqrt} shows that
this assumption was relevant: the solution $v^\e$ associated to $u^\e$
in \eqref{eq:exsqrt} is given by 
\begin{equation*}
v^\e(t,x)= \frac{1}{\e^{n/4}}\V \left(t-t_0,
\frac{x-x_0}{\sqrt\e}\right)\, ,\text{ where } \V =
e^{i\frac{t+t_0}{2}\Delta}\U(-t_0)\, . 
\end{equation*}
For any $T>0$ independent of $\e$, it satisfies \eqref{eq:lincfg} with
$I^\e=[0,T]$, but
$u^\e$ is not linearizable on $[0,T]$ in $L^2$; notice that 
$v^\e$ does not satisfy \eqref{eq:linear}, which is reassuring.

The proof that \eqref{eq:linear} is necessary for linearizability in
$L^2$ relies on profile decomposition for $L^2$ solutions of
\eqref{eq:NLS}. It was established in \cite{MerleVega98} for the case
$n=2$. We prove it in the one-dimensional case in
Section~\ref{sec:profile}. 
\begin{defin}\label{defin:ortho}
If $(h^\e_j,t^\e_j,x^\e_j,\xi^\e_j)_{j \in \N}$ is a  family of
sequences in $\R_+\setminus\{0\}\times\R\times \R^n\times \R^n$, then
we say that 
$(h^\e_j,t^\e_j,x^\e_j,\xi^\e_j)_{j 
\in \N}$  is an  orthogonal family if
\begin{equation*}
\limsup_{\e \rightarrow
0}\left(\frac{h^\e_j}{h^\e_k} +\frac{h^\e_k}{h^\e_j} + \frac{ |t^\e_j -
t^\e_k|}{(h^\e_j)^2}
+ \left| \frac{ x^\e_j -x^\e_k}{h^\e_j} +
\frac{t^\e_j\xi^\e_j -t^\e_k\xi^\e_k }{h^\e_j} 
 \right|
\right) = \infty\, , \quad \forall j \neq k. 
\end{equation*}
\end{defin}
\begin{theo}[Linear profiles]\label{theo:profile}
Let $n=1$ or $2$, and $\U_0^\e$ a bounded family in
$L^2(\R^n)$.\\
\emph{i)} 
Up to extracting a subsequence, there exist an
orthogonal family $(h^\e_j,t^\e_j,x^\e_j,\xi^\e_j)_{j \in \N}$ 
in $ \R_+\setminus\{0\}\times\R\times \R^n\times \R^n$, and a family
$(\phi_j)_{j \in \N}$ bounded in~$L^2(\R^n) $, such that for every
$\ell \ge 1$, 
\begin{equation*}
\begin{aligned}
& e^{i\frac{t}{2}\Delta}\U^\e_0 =\sum_{j=1}^\ell \H_j^\e (\phi_j)(t,x)
+r^\e_\ell(t,x)\, , \\
\text{where }&\ \ \ 
\H_j^\e (\phi_j)(t,x) =e^{i\frac{t}{2}\Delta}\left(e^{ix\cdot
 \xi^\e_j}e^{-i\frac{t^\e_j}{2}\Delta}\frac{1}{(h^\e_j)^{n/2}}\phi_j
 \left(\frac{x-x^\e_j}{h^\e_j} \right) \right)\, ,\\
\text{and }&\ \ \ 
\limsup_{\e \to 0}\|r^\e_\ell\|_{L^{2+4/n}(\R\times\R^n)} \Tend \ell
{+\infty} 0\, . 
\end{aligned}
\end{equation*}
Furthermore, for every $\ell \ge 1$, we have
\begin{equation}\label{eq:orth}
\Norm{\U_0^\e}{L^2(\R)}^2 = \sum_{j=1}^\ell \Norm{\phi_j}{L^2(\R)}^2
+ \Norm{r^\e_\ell}{L^2(\R)}^2 +o(1) \quad \text{as }\e \to 0\, .
\end{equation}
\emph{ii)} If in addition the family $(\U_0^\e)_{0<\e\le 1}$ is
bounded in $H^1(\R^n)$, or more generally if 
\begin{equation}\label{eq:echelle1}
\limsup_{\e\to 0}\int_{|\xi|>R} \left| \widehat \U_0^\e(\xi)\right|^2
d\xi \to 0 \quad \text{as }R\to +\infty\, ,
\end{equation}
then for every $j\ge 1$, $h_j^\e \ge 1$, and $(\xi_j^\e)_\e$ is
bounded, $ |\xi_j^\e|\leq C_j$.
\end{theo}
To state the nonlinear analog to that result, we introduce the
following definition:
\begin{defin}\label{def:profilNL}
Let $\Gamma^\e=(h^\e,t^\e,x^\e,\xi^\e)$ 
be a sequence in  $\R_+\setminus\{0\}\times\R\times \R^n\times \R^n$
such that $t^\e/(h^\e)^2$ has a limit in $[-\infty,+\infty]$ as $\e$ goes to
zero. For $\phi\in 
L^2(\R^n)$, we define the nonlinear profile $\U$ associated to
$(\phi,\Gamma^\e)$ as the unique maximal solution of the nonlinear equation
\eqref{eq:NLS} satisfying
\begin{equation*}
\Norm{\U\(\frac{-t^\e}{(h^\e)^2}\)-
e^{-i\frac{t^\e}{2(h^\e)^2}\Delta}\phi}{L^2(\R^n)}\Tend
\e 0 0\, . 
\end{equation*}
\end{defin}
Essentially, $\phi$ is a Cauchy data for $\U$ if $t^\e/(h^\e)^2$ has a finite
limit, and an asymptotic state (scattering data) otherwise.

\begin{theo}[Nonlinear profiles]\label{theo:profileNL}
Let $n=1$ or $2$, $\U_0^\e$ a bounded family in
$L^2(\R^n)$ and $\U^\e$ the solution to \eqref{eq:NLS} with initial
datum $\U_0^\e$. Let $(\phi_j,\Gamma_j^\e)_{j\in\N^*}$ be the family of
linear profiles 
given by Theorem~\ref{theo:profile}, and $(\U_j)_{j\in\N^*}$ the
family given by Definition~\ref{def:profilNL} (up to
the extraction of a subsequence). \\
Let $I^\e\subset \R$ be a family of open intervals containing the
origin. The following statements are equivalent:
\begin{itemize}
\item[(i)] For every $j\ge 1$, we have
\begin{equation*}
\limsup_{\e \to 0 }\Norm{\U_j}{L^{2+4/n}(I^\e_j\times \R^n)}<+\infty\, ,\quad
\text{where }I_j^\e := (h_j^\e)^{-2}\(I^\e-t_j^\e\)\, . 
\end{equation*}
\item[(ii)] $\dis \limsup_{\e \to 0 }\Norm{\U^\e}{L^{2+4/n}(I^\e\times
\R^n)}<+\infty $.
\end{itemize}
Moreover, if \emph{(i)} or \emph{(ii)} holds, then $\U^\e = \dis
\sum_{j=1}^\ell \U_j^\e + r_\ell^\e + \rho_\ell^\e$, where $r_\ell^\e$ is
given by Theorem~\ref{theo:profile}, and:
\begin{align}
\limsup_{\e \to 0} &\(\Norm{\rho_\ell^\e}{L^{2+4/n}(I^\e\times \R^n)} +
\Norm{\rho_\ell^\e}{L^\infty(I^\e;L^2( \R^n))}\) \Tend \ell {+\infty} 0\,
,\label{eq:resteNL}\\
 \U_j^\e (t,x)&= e^{ix\cdot \xi_j^\e-i\frac{t}{2}(\xi_j^\e)^2}
 \frac{1}{(h_j^\e)^{n/2}} \U_j\( \frac{t-t_j^\e}{(h_j^\e)^2}\virgp
 \frac{x-x_j^\e-t\xi_j^\e}{h_j^\e}\) \label{eq:profNL}.
\end{align}
\end{theo}

We give two applications to these results, besides the proof of
Theorem~\ref{theo:cns}. The first one is the equivalent of 
\cite[Theorem~1.2]{CFG}, which characterizes the obstructions to
linearizability. The second one concerns the properties of blowing up
solutions, in the same spirit as \cite{MerleVega98}. 

The equivalent to \cite[Theorem~1.2]{CFG} is the following. 
\begin{cor}\label{cor:nonlin}
Assume $n=1$ or $2$, and let $u_0^\e$ be bounded in
$L^2(\R^n)$.
Let $T>0$ and assume that \eqref{eq:linear} is not
satisfied with $I^\e=[0,T]$. Then up to the extraction of a
subsequence, there exist an 
orthogonal family $(h^\e_j,t^\e_j,x^\e_j,\xi^\e_j)_{j \in \N}$, 
a family
$(\phi_j)_{j \in \N}$, 
bounded in $L^2(\R^n) $, such that: 
\begin{equation}\label{eq:DACI}
\begin{aligned}
&u^\e_0(x) =\sum_{j=1}^\ell \widetilde H_j^\e (\phi_j)(x)
+w^\e_\ell(x)\, ,&\\
\text{where }&
\widetilde H_j^\e (\phi_j)(x) =e^{ix\cdot
 \xi^\e_j/\sqrt\e}e^{-i\e\frac{t^\e_j}{2}\Delta}\left(
\frac{1}{(h^\e_j\sqrt\e)^{n/2}}\phi_j
 \left(\frac{x-x^\e_j}{h^\e_j\sqrt\e} \right)\right)\, ,&\\
\text{and }&
\limsup_{\e \to 0}\e\|e^{i\e\frac{t}{2}\Delta}
w^\e_\ell\|_{L^{2+4/n}(\R\times\R^n)}^{2+4/n}
\Tend \ell {+\infty} 0\, . &
\end{aligned}
\end{equation}
We have $\liminf  t^\e_j/(h^\e_j)^2 \not = -\infty$,
$\liminf  (T-t^\e_j)/(h^\e_j)^2 \not = -\infty$ (as $\e \to 0$),
and $h_j^\e \leq 1$ for every $j\in \N$. \\
If $t^\e_j/(h^\e_j)^2 \to +\infty$ as $\e \to 0$, then we also have
\begin{equation}\label{eq:DA}
\begin{aligned}
\widetilde H_j^\e (\phi_j)(x) =& e^{i\frac{x\cdot
 \xi^\e_j}{\sqrt\e}+in\frac{\pi}{4}}e^{-i\frac{|x-x^\e_j|^2}{2\e t^\e_j}}\left(
\frac{h^\e_j}{t^\e_j\sqrt\e}\right)^{n/2}\widehat\phi_j
 \left(-\frac{h^\e_j}{t^\e_j\sqrt\e}(x-x^\e_j) \right)\\
&+o(1) \text{ in }L^2(\R^n)\text{ as }\e
    \to 0\, , 
\end{aligned}
\end{equation}
where $\widehat\phi$ stands for the Fourier transform of $\phi$:
$\widehat\phi(\xi) = (2\pi)^{-n/2}\int e^{-ix\cdot
\xi}\phi(x)dx$.\\ 
If in addition $u_0^\e$ is bounded in $H^1_\e$, then
we have $h^\e_j \geq \sqrt\e$.
\end{cor}

\begin{rema*}
Even if $u_0^\e$ is bounded in $H^1_\e$, we cannot say more than $\phi_j\in
L^2(\R^n)$, while in \cite{CFG}, the $H^1_\e$ assumption implied
$\phi_j\in H^1(\R^n)$. This is due to the fact that several scales of
concentrations must be taken into account in the present case, while
in \cite{CFG}, only the scale $\e$ was relevant. In that case, the
profile decomposition in the homogeneous space $\dot H^1(\R^n)$
performed in \cite{Keraani01} could be used to deduce properties in
the \emph{inhomogeneous} Sobolev space $H^1$. In our case, we cannot
compare the $L^2$ and 
$\dot H^1$ profile decompositions.
\end{rema*}
\medskip

\begin{rema*}Compare this result with \cite[Theorem~1.2]{CFG}.\\
\noindent $\bullet$ Scales. As we already
mentioned, not only the scale $\e$ must be considered in the
obstructions to the linearizability in $H^1_\e$, but every scale
between $\e$ and $\sqrt\e$.
Examples \eqref{eq:exsqrt} and
\eqref{eq:excfg} can thus be considered as two borderline
cases.\\
\noindent $\bullet$ Quadratic oscillations. The asymptotic expansion
\eqref{eq:DA} highlights quadratic 
oscillations in the initial data, which are exactly
$\e$-oscillatory, unless $t^\e_j/(h^\e_j)^2$ is
bounded. That case corresponds to initial focusing for $u^\e$ (see for
instance \eqref{eq:exsqrt}). In \cite{CFG}, this phenomenon was
excluded  by the assumption
\begin{equation*}
\e^2 \|u^\e_0\|_{L^{2\si +2}}^{2\si +2} \Tend \e 0 0\, ,
\end{equation*}
because the only relevant concentrating scale was $\e$. In the present case, 
every profile such that $\sqrt\e\ll h^\e_j  \leq 1$ satisfies the above
property, and concentrates with the scale
$h^\e_j\sqrt\e \not =\e$ at time $t=t^\e_j$. It also concentrates with
the same scale at time $t=0$ if $t^\e_j/(h^\e_j)^2$ is
bounded.  So
it is a matter of  choice to consider whether 
or not quadratic oscillations are necessary to have a leading
order nonlinear behavior, according to the way one treats initial
focusing.\\
\noindent $\bullet$ Properties of $t_j^\e$. 
The localization of the cores in time is not as precise as
in \cite{CFG}, where we had $\limsup t_j^\e \in [0,T]$. We actually
have the same condition from the properties $\liminf
t^\e_j/(h^\e_j)^2 \not = -\infty$ and
$\liminf  (T-t^\e_j)/(h^\e_j)^2 \not = -\infty$, provided that the
scale $h^\e_j$ goes to zero as $\e \to 0$. When $h^\e_j$ is constant,
we cannot say much about $t_j^\e$, see \eqref{eq:exsqrt}. \\
\end{rema*}

The second application of Theorem~\ref{theo:profile} concerns finite
time blow up, which may occur for $H^1$-solutions of \eqref{eq:NLS}
when $\lambda =-1$ (not when $\lambda =1$, from the conservation of
energy). For solutions $\U$ in $L^2$ and not necessarily in 
$H^1$, the  
conservation of mass shows that the only obstruction to global
existence in $L^2$ is the unboundness of
$\|\U\|_{L^{2+4/n}([0,T]\times \R^n)}$ (see e.g. \cite{CW89}). 

\begin{cor}\label{cor:blowup}
Assume $n=1$ or $2$. Let $\U$ be an $L^2$-solution to \eqref{eq:NLS},
and assume that $\U$ blows up at time $T>0$ (not 
before),
\begin{equation*}
\int_0^T\!\!\!\! \int_{\R^n}|\U(t,x)|^{2+\frac{4}{n}}dxdt =+\infty\, .
\end{equation*}
Let $(t_k)_{k\in \N}$ be an increasing sequence going to $T$ as $k \to
+\infty$. Then up to a subsequence, there exist $x_j^k,y_j^k\in
\R^n$, $\rho_j^k,h^k_j>0$, $t_j^k\geq 0$ and a family $(\U_j,{\widetilde
\U}_j)_{j\in \N}$ 
bounded in $L^2$ such that
\begin{equation}\label{eq:dvptblowup}
\begin{aligned}
&\U(t_k,x)= \sum_{j=1}^\ell e^{ix\cdot y_j^k}
\frac{1}{(\rho_j^k)^{n/2}}\U_j \left(
\frac{x-x_j^k}{\rho_j^k}\right)\\
& +\sum_{j=1}^\ell e^{ix\cdot y_j^k}
e^{-i\frac{|x-x^k_j|^2}{2(T-t_k)t^k_j}}
\frac{1}{(\widetilde\rho_j^k)^{n/2}}\widetilde\U_j \left( 
\frac{x-x_j^k}{\widetilde\rho_j^k}\right)
 + \W_\ell^k(x)\, ,&\\
\text{with }& \limsup_{k\to +\infty}\left\|
e^{i\frac{t}{2}\Delta}\W_\ell^k\right\|_{L^{2+4/n}(\R\times
\R^n)} \Tend \ell {+\infty} 0\quad ,\quad \widetilde\rho_j^k =
\frac{t^k_j \sqrt{T-t_k}}{h^k_j}\, ,&
\end{aligned}
\end{equation}
and the additional properties, for every $j\in \N$:
\begin{equation}\label{eq:vitesseasym}
\lim_{k\to +\infty} \frac{T-t_k}{(\rho_j^k)^2} \geq 1\, ,
\end{equation}
\begin{equation}\label{eq:profquad}
\text{The sequence }(t_j^k)_{k\in\N}\text{ is bounded, and }\lim_{k\to
+\infty} \frac{t_j^k}{(h^k_j)^2} =+\infty\, .
\end{equation}
Moreover the terms in the sum \eqref{eq:dvptblowup} are pairly orthogonal
in the limit $k\to +\infty$, each term being orthogonal to $\W_\ell^k$.
\end{cor}
\begin{rema*}
For the profiles associated to $\U_j$,  \eqref{eq:vitesseasym} shows 
that the blow up rate is bounded from
below by $(T-t)^{-1/2}$ in the $L^2$ case. In the $H^1$ case, this
property is well known (see \cite{CW90} or \cite{CazCourant}). For the
profiles $\widetilde\U_j$, it is less 
clear. Assume that $\widetilde\U_j$ is smooth, then the $\dot H^1$
norm of the profiles associated to $\widetilde\U_j$ is of order
\begin{equation*}
\left\| \frac{e^{ix\cdot
\xi_j^k-i\frac{|x-x^k_j|^2}{2(T-t_k)t^k_j}}}{(\widetilde\rho_j^k)^{n/2}}
\widetilde\U_j\left(  
\frac{x-x_j^k}{\widetilde\rho_j^k}\right)\right\|_{\dot H^1} \sim
|\xi_j^k| + \frac{1}{h^k_j \sqrt{T-t_k}} + \frac{(h^k_j)^2}{t^k_j}
 \frac{1}{h^k_j \sqrt{T-t_k}}\, .
\end{equation*} 
The second term is due to quadratic oscillations, and dominates the
last term, obtained by differentiating $\widetilde\U_j$, from
\eqref{eq:profquad}. Since from \eqref{eq:profquad}, $h^k_j\to 0$ as
$k\to +\infty$, this suggests that the blow up rate for the profiles
associated to $\widetilde\U_j$ is also bounded from
below by $(T-t)^{-1/2}$ (and is large compared to this minimal rate). 
\end{rema*}

\begin{rema*}
Some blowing up solutions are known explicitly
\cite{Weinstein83}. They are of the form 
\begin{equation*}
\U(t,x)=e^{-i\frac{|x|^2}{2(T-t)} + 
\frac{i}{T-t}}\frac{1}{(T-t)^{n/2}}Q\left(\frac{x}{T-t}\right)\, ,
\end{equation*}
where $Q$ denotes the
unique spherically symmetric solution 
of (see \cite{Strauss77}, \cite{Kwong})
\begin{equation*}
-\frac{1}{2}\Delta Q + Q  =-\lambda|Q|^{4/n}Q\,,  \quad
Q >0 \textrm{ in }\R^n\, .
\end{equation*}
It is proven in \cite{MerleDuke} that up to the invariants of
\eqref{eq:NLS}, these are the only $H^1$ blowing up solutions with
minimal mass $\|\U\|_{L^2}=\|Q\|_{L^2}$. This
yields
\begin{equation*}
\U(T-\e,x) = e^{-i\frac{|x|^2}{2\e} + 
\frac{i}{\e}}\frac{1}{\e^{n/2}}Q\left(\frac{x}{\e}\right)\, ,
\end{equation*}
which is equivalent, up to the extraction of a subsequence, to:
\begin{equation*}
e^{-i\frac{|x|^2}{2\e} + 
i\theta}\frac{1}{\e^{n/2}}Q\left(\frac{x}{\e}\right)\, ,
\end{equation*}
for some $\theta \in \R$. This term may look like a profile
$\widetilde\U_j$, because it contains a quadratic phase, with
$t^\e_j=1$. However, the quadratic oscillation is not
relevant in the profile decomposition. Indeed, 
\begin{equation*}
e^{-i\frac{|x|^2}{2\e} + 
i\theta}\frac{1}{\e^{n/2}}Q\left(\frac{x}{\e}\right)=
e^{ 
i\theta}\frac{1}{\e^{n/2}}Q\left(\frac{x}{\e}\right)+o(1)
\quad\text{in }L^2\, ,
\end{equation*}
and small terms in $L^2$ are linearizable from Strichartz estimates
(see \eqref{eq:strichartzhomogene} below). Although the explicit
formula for these solutions seems to rely on very rigid properties,
our interpretation is as follows. The 
quadratic oscillations gather some mass of $u$ near one point, and
 ignite the blow-up phenomenon: these oscillations, which appear
after a pseudo-conformal transform (see
e.g. \cite{Niederer,GV82,Weinstein86,MerleDuke}), turn a 
non-dispersive solution (typically, a solitary wave) into a
self-focusing solution. A similar explicit formula is available in the
semi-classical limit for \eqref{eq:cfg}, see \cite{Ca2},
Equation~(2.16) and the 
following one; the oscillatory phenomenon seems to be
somehow decoupled from the amplitude one. The quadratic oscillations
correspond to the phase 
given by the geometric optics approach, which describes the geometry
of the propagation. Our point of view is reinforced by this
approach, even though, as mentioned above, explicit formulae may hide
other phenomena (see \cite{BourgainWang,Perelman,MerleRaphaelInv}). 
\end{rema*}

\tableofcontents

\section{Preliminary estimates}
\label{sec:prelim}

First, notice
that the dependence upon $\e$ in \eqref{eq:ck} can be ``removed'' by
the change of unknown function
\begin{equation}\label{eq:chgt}
u^\e(t,x)=\frac{1}{\e^{n/4}}\U^\e \left(t,
\frac{x}{\sqrt\e}\right)\, .
\end{equation}
One checks that $u^\e$ solves \eqref{eq:ck} on $I^\e$ if and only if
$\U^\e$ solves \eqref{eq:NLS} on $I^\e$, and
\begin{equation}\label{eq:equivalence}
\|u^\e(t)\|_{L^2}=\|\U^\e(t)\|_{L^2}\quad ; \quad 
\e\|u^\e(t)\|_{L^{2+4/n}}^{2+4/n}=\|\U^\e(t)\|_{L^{2+4/n}}^{2+4/n}\, .
\end{equation}

In this section, we recall the classical Strichartz estimates, then we
establish 
a refined Strichartz inequality in the space dimension one case. 

\subsection{Classical Strichartz estimates}

The original Strichartz estimate \cite{Strichartz,GV85c}, which holds in any
space dimension, states the following: there
exists a  constant $C$ such that for any $\phi \in L^2(\R^n)$,
\begin{equation}\label{eq:strichartzhomogene}
\left\| e^{i\frac{t}{2}\Delta}\phi \right\|_{L^{2+4/n}(\R \times
\R^n)} \leq C \| \phi\|_{L^2(\R^n)}\, .
\end{equation}
In the case of inhomogeneous Schr\"odinger equations, we have a
similar estimate, which was first proved in \cite{Yajima87}. Denote
\begin{equation*}
\boxed{\gamma = 2+\frac{4}{n}\, ,}
\end{equation*}
and $\gamma '$ its H\"older-conjugate exponent. There
exists a constant $C$ such that for any time interval $I\ni 0$ and any
$\psi \in L^{\gamma '}(I\times \R^n)$,
\begin{equation}\label{eq:strichartzinhomogene}
\begin{aligned}
&\left\|\int_0^t e^{i\frac{t-s}{2}\Delta}\psi(s)ds
\right\|_{L^\gamma(I \times
\R^n)}&\leq C \| \psi \|_{L^{\gamma '}(I\times \R^n)}\, ,\\
&\left\|\int_0^t e^{i\frac{t-s}{2}\Delta}\psi(s)ds
\right\|_{L^\infty(I ;L^2(
\R^n))} &\leq C \| \psi \|_{L^{\gamma '}(I\times \R^n)}\, .
\end{aligned}
\end{equation}
Other estimates are available, but we shall use here only the
three recalled above.

\subsection{A refined Strichartz estimate}

Following \cite{Bourgain95}, a refined Strichartz inequality  was
proved in \cite{MVV} for space dimension two: 
\begin{equation}\label{eq:MVV} 
\begin{aligned}
\left\| e^{i\frac{t}{2}\Delta}u_0\right\|_{L^4(\R_t\times\R^2_x)}
&\lesssim \left\| \widehat {u_0}\right\|_{\chi_p},\quad \text{for }p>
\frac{12}{7}, \quad \text{where}\\
\left\| f\right\|_{\chi_p} &= \( \sum_{j\in \Z} \sum_{\tau\in C_j}
2^{4j} \( \frac{1}{2^{2j}}\int_\tau |f|^p\)^{4/p}\)^{1/4}\, .
\end{aligned}
\end{equation}
Here $\tau$ denotes a square with side length $2^j$, and $C_j$ denotes
a corresponding grid of the plane. This estimate was used
in \cite{Bourgain98,MerleVega98}. We prove its (simpler)
analog in space dimension one.
\begin{prop}\label{prop:refined}
Let $p>1$. There exists $C_p$ such that for any $f\in L^2(\R)$,
\begin{equation*}
\left\| e^{i\frac{t}{2}\d_x^2}f\right\|_{L^6(\R_t\times\R_x)} \le C_p 
\( \sup_{{\tau>0}\atop{\xi_0\in\R}}\tau^{\frac{1}{2}-\frac{1}{p}}
\left\| \widehat
f\right\|_{L^p([\xi_0-\tau,\xi_0+\tau])}\)^{1/3}\|f\|^{2/3}_{L^2(\R)}\,
. 
\end{equation*}
\end{prop}
\begin{proof}
The proof follows very closely the argument used in \cite{KPV00} in
the context of KdV equation. 

By using the explicit formula for the fundamental solution
$e^{i\frac{t}{2}\d_x^2}$, we get
\begin{equation*}
\left|e^{i\frac{t}{2}\d_x^2}f\right|^2 = \iint_{\R^2}
 e^{it(\eta^2-\xi^2)+ix(\xi-\eta)} \widehat f (\xi)\widehat f (\eta)
 d\xi d\eta\, .
\end{equation*}
Introduce the change of variables $u=\eta^2-\xi^2$ and $v=\eta-\xi$:
\begin{equation*}
\left|e^{i\frac{t}{2}\d_x^2}f\right|^2 = \iint_{\R^2}
 e^{itu-ixv} \widehat f (\xi)\widehat f (\eta)
 \frac{dudv}{|\xi-\eta|^{1/2}}\, .
\end{equation*}
We use the usual trick
$\|e^{i\frac{t}{2}\d_x^2}f\|_{L^6(\R^2)}^6 = \|
|e^{i\frac{t}{2}\d_x^2}f|^2 \|_{L^3(\R^2)}^3$. From
Hausdorff--Young's inequality and the inverse change of variable, we
infer
\begin{equation*}
\left\|e^{i\frac{t}{2}\d_x^2}f\right\|_{L^6(\R^2)}^2 \lesssim \(\iint_{\R^2}
 \left|\widehat f (\xi)\right|^{3/2}\left|\widehat f
 (\eta)\right|^{3/2}  \frac{d\xi d\eta}{|\xi-\eta|^{1/2}}\)^{2/3}\, .
\end{equation*}
The end of the proof is analogous to that of \cite[Theorem~3]{KPV00}.
Cauchy--Schwarz inequality yields
\begin{equation*}
\left\|e^{i\frac{t}{2}\d_x^2}f\right\|_{L^6(\R^2)}^3 \lesssim \left\|
\widehat f\right\|_{L^2(\R)} \( 
\int_\R \left|I_{1/2}\( \left| \widehat f
(\xi)\right|^{3/2}\)\right|^2 \left|\widehat f(\xi)\right|d\xi \)^{1/2},
\end{equation*}
where $I_{1/2}$ stands for the fractional integration
$\frac{1}{|x|^{1/2}}\ast$. By Fefferman--Phong's weighted inequality
\cite{Fefferman83}, we get
\begin{equation*}
\( 
\int_\R \left|I_{1/2}\( \left| \widehat f
(\xi)\right|^{3/2}\)\right|^2 \left|\widehat f(\xi)\right|d\xi
\)^{1/2} \le C_p \sup_{{\tau>0}\atop{\xi_0\in\R}}\tau^{\frac{1}{2}-\frac{1}{p}}
\left\| \widehat
f\right\|_{L^p([\xi_0-\tau,\xi_0+\tau])}\|\widehat f\|_{L^2(\R)},
\end{equation*}
which completes the proof of the proposition. 
\end{proof}

\section{Proof of Theorem~\ref{theo:profile}: linear profile decomposition}
\label{sec:profile}

In this section, we prove Theorem~\ref{theo:profile} in the case
$n=1$. The case $n=2$ was established in \cite{MerleVega98}. For the
benefit of the reader, we give a complete proof in the one-dimensional
case.
We follow essentially the same
lines as in \cite{BG3,MerleVega98,Keraani01,IsabelleBSMF}. The idea
relies on an exhaustion 
algorithm inspired from \cite{MetivierSchochet98}, and first used for
such contexts as the present one in \cite{PG98}.

The beginning of the proof does not rely on the assumption $n=1$. We
thus write it with a general $n\ge 1$, and we point out the steps
which are bound to the case $n=1$ (see 
Remark~\ref{rema:dimension}).  
We resume some notations used in the introduction:
\begin{nota}
For a sequence $\Gamma_j^\e=(h^\e_j,t_j^\e,x_j^\e,\xi_j^\e)$ 
in  $\R_+\setminus\{0\}\times\R\times \R^n\times \R^n$, we denote
\begin{equation}\label{eq:defH}
\H_j^\e (\phi_j)(t,x) =e^{i\frac{t}{2}\Delta}\left(e^{ix\cdot
 \xi^\e_j}e^{-i\frac{t^\e_j}{2}\Delta}\frac{1}{(h^\e_j)^{n/2}}\phi_j
 \left(\frac{x-x^\e_j}{h^\e_j} \right) \right).
\end{equation}
For a sequence $\widetilde\Gamma_j^\e=(h^\e_j,x_j^\e,\xi_j^\e)$ 
in  $\R_+\setminus\{0\}\times\R\times \R^n\times \R^n$, we denote
\begin{equation}\label{eq:defHtilde}
\widetilde\H_j^\e (\phi_j)(x) =\H_j^\e (\phi_j)(0,x) =e^{ix\cdot
 \xi^\e_j}e^{-i\frac{t^\e_j}{2}\Delta}\frac{1}{(h^\e_j)^{n/2}}\phi_j
 \left(\frac{x-x^\e_j}{h^\e_j} \right).
\end{equation}
\end{nota}
The following identity is straightforward:
\begin{equation}\label{eq:Hexpl}
\begin{aligned}
\H_j^\e (\phi_j)(t,x) &=e^{ix\cdot \xi_j^\e
-i\frac{t}{2}(\xi^\e_j)^2}  \frac{1}{(h^\e_j)^{n/2}} \V_j \left(
\frac{t-t^\e_j}{(h^\e_j)^{2}}\virgp \frac{x-x^\e_j-t\xi_j^\e}{h^\e_j}\right)
,\\
\text{where }&\V_j(t)=e^{i\frac{t}{2}\Delta} \phi_j\, .
\end{aligned}
\end{equation}
\begin{rema}\label{rem:HH-1}
If two sequences $\Gamma_j^\e=(h^\e_j,t_j^\e,x_j^\e,\xi_j^\e)$ and $
\Gamma_k^\e=(h^\e_k,t_k^\e,x_k^\e,\xi_k^\e)$ are not orthogonal, then, up
to a subsequence, $(\widetilde\H_j^\e)^{-1}\widetilde\H_k^\e \to \widetilde\H$
strongly as $\e\to 0$, where $\widetilde\H$ is isometric on
$L^2(\R^n)$.  
\end{rema}

Let ${\bf \U_0}=(\U_0^\e)_{0<\e\le 1}$ be a bounded sequence in
$L^2(\R^n)$. We denote by $\mathcal V ({\bf \U_0})$ the set of
weak limits of subsequences of the
form $(\widetilde\H^\e)^{-1}\U_0^\e$ for some $\Gamma^\e =
(h^\e,t^\e,x^\e,\xi^\e)$ 
in $\R_+\setminus\{0\}\times\R\times \R^n\times \R^n$:
\begin{equation*}
\mathcal V ({\bf \U_0}) = \left\{ w-\lim_{k\to +\infty 
}(\widetilde\H^{\e_k})^{-1}\U_0^{\e_k}\ ; \ \e_k\Tend
k{+\infty} 0,\  \Gamma^{\e_k}\in
\R_+\setminus\{0\}\times\R\times \R^n\times \R^n \right\}.
\end{equation*}
We denote
\begin{equation*}
\eta({\bf \U_0}) =\sup \left\{ \|\phi\|_{L^2(\R^n)}\quad ;\quad \phi
\in \mathcal V ({\bf \U_0})\right\}.
\end{equation*}
We have obviously\footnote{We have the $\limsup$ -- not the
$\liminf$ -- because 
we consider all possible subsequences $\e_k$.}
\begin{equation*}
\eta({\bf \U_0}) \le \limsup_{\e \to 0}\Norm{\U_0^\e}{L^2(\R^n)}\, .
\end{equation*}
We prove that there exist a sequence $(\phi_j)_{j\ge 1}$ and a family
of pairwise orthogonal sequences $\Gamma_j^\e
=(h^\e_j,t_j^\e,x_j^\e,\xi_j^\e)$ such that, up to extracting a
subsequence:
\begin{equation}\label{eq:decomp1}
\U_0^\e = \sum_{j=1}^\ell \widetilde \H_j^\e(\phi_j) + \V_\ell^\e\,
,\quad \text{with }\eta({\bf \V_\ell})\Tend \ell {+\infty} 0\,,
\end{equation}
and with the almost orthogonality identity:
\begin{equation}\label{eq:decomp1'}
\Norm{\U_0^\e}{L^2(\R^n)}^2 =
\sum_{j=1}^\ell\Norm{\phi_j}{L^2(\R^n)}^2 +
\Norm{\V_\ell^\e}{L^2(\R^n)}^2 +o(1)\quad \text{as }\e \to 0\, .
\end{equation}
Indeed, if $\eta({\bf \U_0})=0$, then we can take $\phi_j\equiv 0$ for
all $j$. Otherwise, we choose $\phi_1\in \mathcal V ({\bf \U_0})$ such
that 
\begin{equation*}
\Norm{\phi_1}{L^2(\R^n)}\ge \frac{1}{2}\eta({\bf \U_0})>0\,.
\end{equation*}
By definition, there exists some sequence $\Gamma_1^\e =
(h^\e_1,t_1^\e,x_1^\e,\xi_1^\e)$ such that, up to extracting a
subsequence, we have:
\begin{equation*}
\( \widetilde \H_1^\e\)^{-1}\U_0^\e \rightharpoonup \phi_1\, .
\end{equation*}
We set $\V_1^\e = \U_0^\e - \widetilde \H_1^\e(\phi_1)$, and we get
\begin{equation*}
\Norm{\U_0^\e}{L^2(\R^n)}^2 =
\Norm{\phi_1}{L^2(\R^n)}^2 +
\Norm{\V_1^\e}{L^2(\R^n)}^2 +o(1)\quad \text{as }\e \to 0\, .
\end{equation*}
Now, we replace $\U_0^\e$ with $\V_1^\e$, and repeat the same
process. If $\eta ({\bf \V_1})>0$, we get $\phi_2$, $\Gamma_2^\e$ and
${\bf \V_2}$. Moreover, $\Gamma_1^\e$ and $\Gamma_2^\e$ are
orthogonal. Otherwise, up to extracting a subsequence, we use
Remark~\ref{rem:HH-1}:
$(\widetilde\H_2^\e)^{-1}\widetilde\H_1^\e \to \widetilde\H$ 
strongly as $\e\to 0$, where $\widetilde\H$ is isometric on
$L^2(\R^n)$. Since
\begin{equation*}
(\widetilde\H_2^\e)^{-1}\V_1^\e =
\((\widetilde\H_2^\e)^{-1}\widetilde\H_1^\e\) (\widetilde\H_1^\e)^{-1}
\V_1^\e  
\end{equation*}
and $(\widetilde\H_1^\e)^{-1}
\V_1^\e  $ converges weakly to zero, this implies $\phi_2\equiv 0$, hence
$\eta({\bf \V_1})=0$, which yields a contradiction. 

Iterating this argument, a diagonal process yields a family of
pairwise orthogonal sequences $\Gamma_j^\e$, and $(\phi_j)_{j\ge 1}$
satisfying \eqref{eq:decomp1'}. Since
$(\U_0^\e)_{0<\e\le 1}$ is bounded in $L^2(\R^n)$, \eqref{eq:decomp1'}
yields
\begin{equation*}
\sum_{j=1}^\ell\Norm{\phi_j}{L^2(\R^n)}^2 \le \limsup_{\e \to 0}
\|\U_0^\e\|_{L^2(\R^n)}^2 \, .
\end{equation*}
Since the bound is independent of $\ell\ge 1$, the series
$\sum\Norm{\phi_j}{L^2(\R^n)}^2  $ is convergent, and 
\begin{equation*}
\Norm{\phi_j}{L^2(\R^n)}\to 0 \quad \text{as } j\to +\infty\, .
\end{equation*}
Furthermore, we have by construction
\begin{equation*}
\eta({\bf \V_\ell}) \leq \Norm{\phi_{\ell-1}}{L^2(\R^n)}\, ,
\end{equation*}
which yields \eqref{eq:decomp1}. \\

When the initial data satisfy \eqref{eq:echelle1}, we alter the above
algorithm. We \emph{impose} the lower bound on the scales and the
boundedness of the cores in the Fourier side:
\begin{align*}
\widetilde{\mathcal V} ({\bf \U_0}) = \Big\{ &w-\lim_{k\to +\infty 
}(\widetilde\H^{\e_k})^{-1}\U_0^{\e_k}\ ; \ \e_k\Tend
k{+\infty} 0,\  \Gamma^{\e_k}\in
[1,+\infty[\times\R\times \R^n\times \R^n\ ,\\
& \text{with
}|\xi^\e|\lesssim 1\Big\}.
\end{align*}
Notice that the assumption $h^\e\ge 1$ is nothing but boundedness away
from zero. Up to an $\e$-independent dilation of the profiles, we may
always assume that a scale bounded away from zero is bounded from
below by $1$. 

Repeating the same algorithm as above, the property
\eqref{eq:echelle1} remains at each step. Notice that the stronger
assumption $\U_0^\e\in H^1$ is not stable: we may not have $\phi_1\in
H^1$. The same lines yield \eqref{eq:decomp1}, and
\eqref{eq:decomp1'}, with $\eta$ replaced by $\widetilde \eta$ with the
natural definition for $\widetilde \eta$. In particular, $\V_\ell^\e$
satisfies \eqref{eq:echelle1} for any $\ell \ge 1$.

Theorem~\ref{theo:profile} stems from
the following proposition:
\begin{prop}\label{prop:2}
We assume $n=1$. Let $(\U_0^\e)_{0<\e\le 1}$ be a family of $L^2(\R)$
such that
\begin{equation*}
\Norm{\U_0^\e}{L^2(\R)}\leq M\quad \text{and}\quad
\Norm{e^{i\frac{t}{2}\d_x^2}\U_0^\e}{L^6(\R_t\times \R_x)}\ge m>0\, . 
\end{equation*}
There exists $\Gamma^\e=(h^\e,t^\e,x^\e,\xi^\e)$ such that, up to a
subsequence,  
\begin{equation*}
\(\widetilde \H^\e\)^{-1}(\U_0^\e)\rightharpoonup \phi\, ,\quad
\text{where } \Norm{e^{i\frac{t}{2}\d_x^2}\phi}{L^6(\R_t\times
\R_x)}\ge \beta(m)>0\, .
\end{equation*}
Moreover, if $(\U_0^\e)_{0<\e\le 1}$ satisfies \eqref{eq:echelle1}, then one
can choose $h^\e \ge 1$ and $|\xi^\e| \lesssim 1$.
\end{prop}
\begin{rema*}
The dependence of $\beta$ upon $M$ is not mentioned in the above
statement. Simply recall that from Strichartz inequality
\eqref{eq:strichartzhomogene}, $m \lesssim M$. 
\end{rema*}
This proposition, together with \eqref{eq:decomp1}, yields
Theorem~\ref{theo:profile}. Indeed, if $\V_\ell^\e = r^\e_\ell(0,x)$
was  such that 
\begin{equation*}
\limsup_{\e \to 0}\Norm{e^{i\frac{t}{2}\d_x^2}\V_\ell^\e}{L^6(\R^2)}
=  \limsup_{\e \to 0}\Norm{r_\ell^\e}{L^6(\R^2)}\not\to 0\quad\text{as
}\ell \to{+\infty} \, ,
\end{equation*}
then there would exist $\ell_k\to +\infty $ as $k\to +\infty$, and $m>0$,
such that for any $k\in \N$,
\begin{equation*}
\limsup_{\e \to 0}\Norm{e^{i\frac{t}{2}\d_x^2}\V_{\ell_k}^\e}{L^6(\R^2)}
\ge m\, .
\end{equation*}
From \eqref{eq:decomp1'}, we have
\begin{equation*}
\limsup_{\e \to 0}\Norm{\V_{\ell_k}^\e}{L^2(\R)}
\le \limsup_{\e \to 0}\Norm{\U_0^\e}{L^2(\R)}=:M\, .
\end{equation*}
From Proposition~\ref{prop:2}, there exists
$\varphi_\ell\in
\eta({\bf \V_\ell})$ with 
\begin{equation*}
\Norm{e^{i\frac{t}{2}\Delta}\varphi_{\ell_k}}{L^6(\R_t\times
\R_x)}\ge \beta(m)>0\, .
\end{equation*}
This implies $\eta({\bf \V_{\ell_k}})\ge \beta(m)>0$ for any $k\in\N$,
which contradicts \eqref{eq:decomp1}.
\begin{proof}
The proof of Proposition~\ref{prop:2} relies on several intermediary
results. First, we extract scales $h^\e_j$ and cores on the Fourier
side $\xi^\e_j$ and obtain a remainder arbitrarily small thanks to the
refined Strichartz estimate.
\begin{lem}\label{lem:1}
Let $(\U_0^\e)_{\e}$ be a bounded sequence in $L^2(\R)$. Then
for every $\delta>0$, there exist $N=N(\delta)$, a family
$(h^\e_j,\theta^\e_j)_{1\le j\le N}\in \R_+\setminus\{0\}\times\R$, and a
family $({\bf g}_j)_{1\le j\le N}$ of bounded sequences in $L^2(\R)$
such that, up to a subsequence,
\begin{itemize}
\item[(i)] If $j \neq k$, $\dis 
\frac{h^\e_j}{h^\e_k} +\frac{h^\e_k}{h^\e_j} + \left| \theta_k^\e -
\frac{h_k^\e}{h_j^\e} \theta_j^\e \right| \Tend \e 0 +\infty$.
\item[(ii)] For every $1\le j\le N$, there exists $F_j$ bounded, compactly
supported, such that
\begin{equation}\label{eq:10}
\sqrt{h_j^\e}\left|  \widehat{g}_j^\e\(h_j^\e \xi+\theta_j^\e\) \right|\le
F_j(\xi)\, . \
\end{equation}
\item[(iii)] For every $\ell \ge 1$ and $x\in\R$, 
\begin{equation}\label{eq:11}
\U_0^\e =
\sum_{j=1}^N g_j^\e + q^\e\, ,\quad \text{with}\quad
\Norm{e^{i\frac{t}{2}\d_x^2}q^\e}{L^6(\R^2)}\le \delta \, .
\end{equation}
\end{itemize}
Moreover, we have the almost orthogonality identity:
\begin{equation}\label{eq:12}
\Norm{\U_0^\e}{L^2}^2 = \sum_{j=1}^N\Norm{g_j^\e}{L^2}^2 +
\Norm{q^\e}{L^2}^2 +o(1)\quad \text{as }\e \to 0\, .
\end{equation}
\end{lem}
\begin{proof}[Proof of Lemma~\ref{lem:1}]
For $\gamma^\e = (h^\e,\theta^\e)\in \R_+\setminus\{0\}\times\R$, we
define
\begin{equation*}
\G^\e (f)(\xi)=\sqrt{h^\e}  f^\e\(h^\e \xi+\theta^\e\) \, .
\end{equation*}
If   $\Norm{e^{i\frac{t}{2}\d_x^2}\U_0^\e}{L^6(\R^2)}\le \delta $,
then nothing is to be proved. Otherwise, up to extracting a subsequence, 
$\Norm{e^{i\frac{t}{2}\d_x^2}\U_0^\e}{L^6(\R^2)}> \delta $. Apply
Proposition~\ref{prop:refined} with $p=4/3$; there exists a family of
intervals $I^\e = [\theta^\e-h^\e,\theta^\e+h^\e]$ such that
\begin{equation*}
\int_{I^\e} \left| \widehat{\U}_0^\e\right|^{4/3} \ge C
\delta^4 (h^\e)^{1/3}\, ,
\end{equation*}
where the constant $C$ is uniform since $(\U_0^\e)_\e$ is bounded in
$L^2$. For any $A>0$, we have
\begin{equation*}
\int_{I^\e\cap \{|\widehat{\U}_0^\e|>A\}} \left|
\widehat{\U}_0^\e\right|^{4/3}\le A^{-2/3}\Norm{ \widehat{\U}_0^\e
}{L^2}^2\, .
\end{equation*}
Taking $A= C'/ \( \sqrt{h^\e}\delta^6\)$ yields
\begin{equation*}
\int_{I^\e\cap \{|\widehat{\U}_0^\e|\le A\}} \left|
\widehat{\U}_0^\e\right|^{4/3}\gtrsim \delta^4 (h^\e)^{1/3}\, . 
\end{equation*}
From H\"older's inequality, we infer
\begin{equation*}
\int_{I^\e\cap \{|\widehat{\U}_0^\e|\le A\}} \left|
\widehat{\U}_0^\e\right|^{2}\ge C''  \delta^6 \, , 
\end{equation*}
for some uniform constant $C''$. Define $\v_1^\e$ and $\gamma_1^\e$ by 
\begin{equation*}
\widehat{\v}_1^\e = \widehat{\U}_0^\e \1_{I^\e \cap 
\{|\widehat{\U}_0^\e|\le A\}} \quad ;\quad \gamma_1^\e =
(h^\e,\theta^\e)\, .
\end{equation*}
We have
\begin{equation*}
\left| \G_1^\e (\widehat{\v}_1^\e )(\xi)\right| \le C(\delta)
\1_{[-1,1]}(\xi)\, ,
\end{equation*}
which is \eqref{eq:10} with $g_j^\e$ replaced by ${\v}_1^\e$. Furthermore,
\begin{equation*}
\Norm{\U_0^\e}{L^2}^2 = \Norm{\U_0^\e- \v_1^\e}{L^2}^2 + \Norm{
\v_1^\e}{L^2}^2\, ,
\end{equation*}
since the supports are disjoint from the Fourier side. 

We repeat the same argument with $\U_0^\e- \v_1^\e$ in place of
$\U_0^\e$. At each step, the $L^2$ norm decreases of at least
$(C'')^{1/2}\delta^3$, with the same constant $C''$ as for the first
step. After $N(\delta)$ steps, we obtain
$(\v_j^\e)_{1\le j\le N(\delta)}$ and $(\gamma_j^\e)_{1\le j\le
N(\delta)}$ satisfying \eqref{eq:10}, such that
\begin{equation}\label{eq:21}
\U_0^\e = \sum_{j=1}^{N(\delta)} \v_j^\e +q^\e\, , \quad \text{with }
\Norm{e^{i\frac{t}{2}\d_x^2}q^\e}{L^6(\R^2)}\le \delta \, ,
\end{equation}
and $\Norm{\U_0^\e}{L^2}^2 =
\sum_{j=1}^{N(\delta)}\Norm{\v_j^\e}{L^2}^2 + \Norm{q^\e}{L^2}^2 +o(1)$
as $\e \to 0$. However, the properties of the first point of the lemma
need not be satisfied. To obtain these properties, we reorganize the
decomposition. 
We say that $\gamma_j^\e$ and $\gamma_k^\e$ are orthogonal if 
\begin{equation*}
\frac{h_j^\e}{h_k^\e} + \frac{h_k^\e}{h_j^\e}+\left| \theta_k^\e -
\frac{h_k^\e}{h_j^\e} \theta_j^\e \right| \to +\infty \quad \text{as }\e
\to 0\, .
\end{equation*}
Define 
\begin{equation*}
g_1^\e  = \sum_{j=1}^{N(\delta)}\v_j^\e - \sum_{\gamma_j^\e \perp
\gamma_1^\e}\v_j^\e\, .
\end{equation*}
If there exists $2\le j_0\le N(\delta)$ such that $\gamma_{j_0}^\e$ is
orthogonal to $\gamma_1^\e$, then we define
\begin{equation*}
g_2^\e  = \sum_{j=1}^{N(\delta)}\v_j^\e - \sum_{{\gamma_j^\e \perp
\gamma_1^\e}\atop {\gamma_j^\e \perp
\gamma_{j_0}^\e}}\v_j^\e\, .
\end{equation*}
Repeating this argument a finite number of times, we rearrange the
above sum. The almost orthogonality relation \eqref{eq:12} holds, since the 
supports of the functions we consider are disjoint from the Fourier
side. Finally, we must make sure that up to an extraction, the first
point of the lemma is satisfied, and that \eqref{eq:10} holds. 

The $\v_j^\e$'s kept in the definition of $g_1^\e$ are such that the
$\gamma_j^\e$ are not orthogonal one to another. It is sufficient to
show that up to an extraction, $\G_1^\e(\v_j^\e)$ is bounded by a
compactly supported bounded function, for such $j$'s. By construction,
$\G_j^\e(\v_j^\e)$ is bounded by a 
compactly supported bounded function; we have
\begin{equation*}
\G_1^\e (\G_j^\e)^{-1} f(\xi) = \sqrt{\frac{h_1^\e}{h_j^\e}}f \(
\frac{h_1^\e}{h_j^\e} \xi +\theta_1^\e -
\frac{h_1^\e}{h_j^\e}\theta_j^\e\).
\end{equation*}
Since $\gamma_j^\e \not\perp \gamma_1^\e$, up to an extraction,
$h_1^\e/ h_j^\e\to \lambda_{1j}\in\R_+\setminus\{0\}$ and $\theta_1^\e -
\frac{h_1^\e}{h_j^\e}\theta_j^\e$ is bounded as $\e\to 0$, which yields
the desired estimate for $\G_1^\e(\v_j^\e)$. Reasoning the same way
for the other terms proves (i) and (ii), and completes the proof of
the lemma.
\end{proof}
Next, we study sequences whose scale $h^\e$ is
fixed, equal to $1$, and extract cores in space-time. 
\begin{prop}\label{prop:3}
Let ${\bf P}=(P^\e)_{0<\e\le 1}$ be a sequence 
such that 
\begin{equation}\label{eq:aspect}
\left|\widehat{P^\e}(\xi)\right| \leq F(\xi)\, ,
\end{equation}
where $F\in L^\infty(\R)$ is compactly supported. Then there exist a 
subsequence of $P^\e$ (still denoted $P^\e$), a family $({\bf
x}_\alpha,{\bf s}_\alpha)_{\alpha\ge 1}$ of sequences in $\R\times\R$,
and a sequence $({\bf \phi}_\alpha)_{\alpha\ge 1}$ of $L^2$ functions,
such that:
\begin{itemize}
\item[(i)] If $\alpha \neq \beta$, $|x_\alpha^\e - x_\beta^\e| +
|s_\alpha^\e - s_\beta^\e| \to +\infty$ as $\e \to 0$.
\item[(ii)] For every $A\ge 1$ and every $x\in\R$, we have:
\begin{equation*}
P^\e(x) = \sum_{\alpha =1}^A e^{-is_\alpha^\e\d_x^2
}\phi_\alpha (x-x_\alpha^\e) + P_A^\e(x)\, ,\quad \text{with}
\end{equation*}
\begin{equation}\label{eq:13}
\limsup_{\e \to 0}\Norm{e^{i\frac{t}{2}\d_x^2}P_A^\e}{L^6(\R^2)} \Tend
A \infty 0\, ,\quad \text{and}
\end{equation}
\begin{equation}\label{eq:14}
\Norm{P^\e}{L^2}^2 = \sum_{\alpha =1}^A \Norm{\phi_\alpha}{L^2}^2 + 
\Norm{P_A^\e}{L^2}^2+o(1)\quad \text{as }\e \to 0\, .
\end{equation}
\end{itemize}
\end{prop}
\begin{proof}[Proof of Proposition~\ref{prop:3}]
Let $\mathcal W ({\bf P})$ be the set of weak limits of subsequences
of ${\bf P}$ after translation in the phase space:
\begin{equation*}
\mathcal W ({\bf P}) = \left\{ w-\lim_{k \to +\infty} e^{i
s^{\e_k}\d_x^2} P^\e (\cdot + x^{\e_k})\ ;\ \e_k\Tend k {+\infty} 0,
\ (x^\e,s^\e)\in \R\times\R\right\}.
\end{equation*}
We denote
\begin{equation*}
\mu ({\bf P}) = \sup\left\{ \Norm{\phi}{L^2}\ ;\ \phi \in
\mathcal W ({\bf P}) \right\}.
\end{equation*}
As in the beginning of this section, we have
\begin{equation*}
\mu ({\bf P}) \le \limsup_{\e \to 0}\Norm{P^\e}{L^2}\, ,
\end{equation*}
and, up to extracting a subsequence, we can write
\begin{equation*}
P^\e(x) = \sum_{\alpha =1}^A e^{-is_\alpha^\e\d_x^2
}\phi_\alpha (x-x_\alpha^\e) + P_A^\e(x)\, ,\quad \mu ({\bf P}_A)\Tend
A {+\infty} 0\, ,
\end{equation*}
with the almost orthogonality identity \eqref{eq:14}. To complete the
proof of Proposition~\ref{prop:3}, we have to prove \eqref{eq:13}. 

Notice that the orthogonality argument yields a result more precise
than \eqref{eq:14}: for every $\alpha \ge 1$ and every $\psi\in
\F(C_0^\infty(\R))$, 
\begin{equation*}
\Norm{\widehat\psi\widehat P^\e}{L^2}^2 = \sum_{\alpha =1}^A
\Norm{\widehat\psi\widehat \phi_\alpha}{L^2}^2 +  
\Norm{\widehat\psi\widehat P_A^\e}{L^2}^2+o(1)\quad \text{as }\e \to 0\, .
\end{equation*}
This fact, together with the assumption \eqref{eq:aspect}, proves that
for every $A\ge 1$, $\widehat P_A^\e$ is supported in $\supp F$, and
\begin{equation}\label{eq:16}
\limsup_{\e \to 0}\Norm{\widehat P_A^\e}{L^\infty}\le \Norm{\widehat
F}{L^\infty} \, .
\end{equation}
Introduce a cut-off
 $\chi(t,x)=\chi_1(t)\chi_2(x)$, with $\chi_j\in\S(\R)$, such that:
\begin{equation*}
\left| \widehat \chi_1\right| + \left| \widehat \chi_2\right| \le 2\quad
 ;\quad \widehat \chi_2\equiv 1\ \text{ on }\supp F\quad
 ;\quad \widehat \chi_1\( \frac{-\xi^2}{2}\) \equiv 1 \
 \text{ on }\supp \widehat \chi_2\, . 
\end{equation*}
Let $\ast$ denote the convolution in $(t,x)$, and $\psi_A^\e(t,x)=
e^{i\frac{t}{2}\d_x^2}P_A^\e$. The function $\chi \ast 
\psi_A^\e$ solves the 
linear Schr\"odinger equation, so
\begin{equation*}
\F \( \chi\ast \psi_A^\e\big|_{t=0}\)(\xi) =
\widehat \chi_1\( 
\frac{-\xi^2}{2}\) \widehat
\chi_2(\xi)\widehat{P_A^\e}(\xi)=\widehat{P_A^\e}(\xi)\, ,
\end{equation*}
from the assumptions on $\chi_1$ and $\chi_2$. Therefore, $\chi\ast \psi_A^\e
=\psi_A^\e$. We use a restriction result in space dimension $1$ (see
e.g. \cite{TaoUtah03}): for every $4<q<6$ and every $\widehat G\in L^\infty(B(0,R))$, 
\begin{equation}\label{eq:restriction}
\Norm{\int_{B(0,R)}e^{i\frac{t}{2}|\xi|^2+ix\cdot
\xi}\widehat G(\xi)d\xi}{L^q(\R^2)} \le C(q,R) \Norm{\widehat G}{L^\infty}\, .
\end{equation}
Fix $4<q<6$. Using \eqref{eq:16} and \eqref{eq:restriction}, we have
\begin{align*}
\limsup_{\e\to 0}\Norm{\chi\ast \psi_A^\e}{L^6(\R^2)} &\le
\limsup_{\e\to 0}\Norm{\chi\ast
\psi_A^\e}{L^q(\R^2)}^{q/6}\limsup_{\e\to 0}\Norm{\chi\ast
\psi_A^\e}{L^\infty(\R^2)}^{1-q/6} \\
&\le \Norm{F}{L^\infty(\R)}^{q/6}\limsup_{\e\to 0}\Norm{\chi\ast
\psi_A^\e}{L^\infty(\R^2)}^{1-q/6}\, .
\end{align*}
On the other hand, the definition of $\mathcal W ({\bf P}_A)$ implies
\begin{equation*}
\limsup_{\e\to 0}\Norm{\chi\ast
\psi_A^\e}{L^\infty(\R^2)}\le \sup\left\{ \left| \iint
\chi(-t,-x)e^{i\frac{t}{2}\d_x^2}\phi \, dxdt\right| \ ; \ \phi \in
\mathcal W ({\bf P}_A)\right\}.
\end{equation*}
Using H\"older's inequality, then Strichartz estimate, we obtain
\begin{align*}
\limsup_{\e\to 0}\Norm{\chi\ast
\psi_A^\e}{L^\infty(\R^2)} &\le \Norm{\chi}{L^{6/5}(\R^2)}\sup\left\{ 
\Norm{e^{i\frac{t}{2}\d_x^2}\phi}{L^6(\R^2)}\ ; \ \phi \in
\mathcal W ({\bf P}_A)\right\}\\
&\lesssim \Norm{\chi}{L^{6/5}(\R^2)} \mu ({\bf P}_A)\, .
\end{align*}
Therefore, 
\begin{equation*}
\limsup_{\e\to 0}\Norm{e^{i\frac{t}{2}\d_x^2}P_A^\e}{L^6(\R^2)}
\lesssim \mu ({\bf P}_A)^{1-\frac{q}{6}}\to 0 \quad \text{as }A\to
+\infty\, ,
\end{equation*}
which completes the proof of Proposition~\ref{prop:3}.
\end{proof}
We can now finish the proof of Proposition~\ref{prop:2}. Back to the
decomposition \eqref{eq:11}, we set, for $1\le j\le N$, 
\begin{equation*}
P_j^\e(x) = e^{-ix \theta_j^\e/h_j^\e }\sqrt{h_j^\e}
g_j^\e\(h_j^\e x \) . 
\end{equation*}
Since $g_j^\e$ satisfies \eqref{eq:10}, the sequence
$(P_j^\e)_{0<\e\le 1}$ satisfies the assumptions of
Proposition~\ref{prop:3}. Thus, for every $1\le j\le N$, there exists
a family $(\phi_{j,\alpha})_{\alpha\ge 1}$ of $L^2$ functions, and a
family $(y_{j,\alpha}^\e,s_{j,\alpha}^\e)\in\R\times\R$, such that
\begin{equation}\label{eq:18}
P_j^\e(x) = \sum_{\alpha =1}^A e^{-is_{j,\alpha}^\e\d_x^2
}\phi_{j,\alpha} \(x-y_{j,\alpha}^\e\) + P_{j,A}^\e(x)\, ,
\end{equation}
together with \eqref{eq:13} and \eqref{eq:14}. For each $1\le j\le N$,
choose $A_j$ such that for $A\ge A_j$, 
\begin{equation*}
\limsup_{\e \to 0} \Norm{e^{i\frac{t}{2}\d_x^2}P_{j,A}^\e}{L^6(\R^2)}
\le \frac{\delta}{N}\, .
\end{equation*}
In terms of $g_j^\e$, \eqref{eq:18} reads
\begin{equation*}
\begin{aligned}
&g_j^\e = \sum_{\alpha =1}^A \widetilde\H^\e_{j,\alpha}
\(\phi_{j,\alpha}\) + \w_{j,A}^\e\, ,\quad \text{where}\\
\Gamma^\e_{j,\alpha} = & \( h_{j}^\e, 2 s_{j,\alpha}^\e,
h_{j}^\e y_{j,\alpha}^\e, \frac{\theta_{j}^\e}{h_j^\e}\)\ ; \quad
\w_{j,A}^\e(x) =
\frac{e^{-ix\xi_{j}^\e}}{\sqrt{h_j^\e}}P_{j,A}^\e \(
\frac{x}{h_j^\e}\).
\end{aligned} 
\end{equation*}
Using \eqref{eq:11}, it follows that 
\begin{equation*}
\U_0^\e = \sum_{j=1}^N\( \sum_{\alpha =1}^{A_j} \widetilde\H^\e_{j,\alpha}
\(\phi_{j,\alpha}\) + \w_{j,A}^\e\) + q^\e\, .
\end{equation*}
Relabeling the pairs $(j,\alpha)$, we get
\begin{equation*}
\U_0^\e = \sum_{j=1}^K \widetilde\H^\e_j \(\phi_j\) + \W^\e\, ,
\end{equation*}
where $K= \dis \sum_{j=1}^N A_j$ and $\dis \W^\e = \sum_{j=1}^N
\w_{j,A_j}^\e +q^\e$. The remainder satisfies 
\begin{equation*}
\limsup_{\e \to 0}\Norm{e^{i\frac{t}{2}\d_x^2}\W^\e}{L^6(\R^2)}\le
2\delta\, .
\end{equation*}
It is clear that the $\Gamma_j^\e$'s are pairwise
orthogonal. Combining \eqref{eq:12} and \eqref{eq:14}, we obtain
\begin{equation*}
\Norm{\U_0^\e}{L^2}^2 = \sum_{j=1}^N \(\sum_{\alpha =1}^{A_j}
 \Norm{\phi_{j,\alpha}}{L^2}^2 +\Norm{\w_{j,A}^\e}{L^2}^2\) +
 \Norm{q^\e}{L^2}^2 +o(1)\quad \text{as }\e \to 0\, .
\end{equation*}
Thus,
\begin{equation}\label{eq:19}
\sum_{j=1}^K \
 \Norm{\phi_{j}}{L^2}^2 \le \limsup_{\e \to 0}\Norm{\U_0^\e}{L^2}^2\le
 M^2\, .
\end{equation}
Since $\dis \Norm{e^{i\frac{t}{2}\d_x^2}\U_0^\e}{L^6(\R^2)}\ge m>0$,
choose $\delta$ small enough  so that
\begin{equation*}
\frac{1}{2}\Norm{e^{i\frac{t}{2}\d_x^2}\U_0^\e}{L^6(\R^2)}^6 \le
\Norm{\sum_{j=1}^K \H_j^\e \(\phi_{j}\)}{L^6(\R^2)}^6\le
\Norm{e^{i\frac{t}{2}\d_x^2}\U_0^\e}{L^6(\R^2)}^6  \, .
\end{equation*}
A classical argument of orthogonality (see e.g. \cite{PG98}) yields,
as $\e \to 0$, 
\begin{equation*}
\Norm{\sum_{j=1}^K \H_j^\e \(\phi_{j}\)}{L^6(\R^2)}^6 =
\sum_{j=1}^K\Norm{ \H_j^\e \(\phi_{j}\)}{L^6(\R^2)}^6 +o(1)=
\sum_{j=1}^K\Norm{ e^{i\frac{t}{2}\d_x^2}\phi_{j}}{L^6(\R^2)}^6 +o(1)\, .  
\end{equation*}
Let $j_0$ be such that $\Norm{
e^{i\frac{t}{2}\d_x^2}\phi_{j_0}}{L^6(\R^2)}= \dis\max_{1\le j\le K}
\Norm{ e^{i\frac{t}{2}\d_x^2}\phi_{j}}{L^6(\R^2)}$. Using Strichartz
estimate, we infer
\begin{align*}
\frac{m^6}{2} &\le \sum_{j=1}^K\Norm{
e^{i\frac{t}{2}\d_x^2}\phi_{j}}{L^6(\R^2)}^6  \le 
\Norm{e^{i\frac{t}{2}\d_x^2}\phi_{j_0}}{L^6(\R^2)}^4 \sum_{j=1}^K\Norm{
e^{i\frac{t}{2}\d_x^2}\phi_{j}}{L^6(\R^2)}^2\\
&\lesssim \Norm{e^{i\frac{t}{2}\d_x^2}\phi_{j_0}}{L^6(\R^2)}^4
\sum_{j=1}^K\Norm{ 
\phi_{j}}{L^2(\R^2)}^2 \lesssim M^2
\Norm{e^{i\frac{t}{2}\d_x^2}\phi_{j_0}}{L^6(\R^2)}^4 \,,
\end{align*}
where the last estimate follows from \eqref{eq:19}. Thus,
\begin{equation*}
\Norm{e^{i\frac{t}{2}\d_x^2}\phi_{j_0}}{L^6(\R^2)} \ge \beta\thickapprox
\frac{m^{3/2}}{M}\, \cdot
\end{equation*}
The pairwise orthogonality of the $\Gamma_j^\e$'s yields
\begin{equation*}
\(\widetilde \H_{j_0}^\e\)^{-1} \U_0^\e \rightharpoonup \phi =
\phi_{j_0} + \W\, ,
\end{equation*}
where $\W$ is the weak limit of $(\widetilde \H_{j_0}^\e)^{-1}
\W^\e$. Since 
\begin{equation*}
\Norm{e^{i\frac{t}{2}\d_x^2}\W}{L^6(\R^2)}\le \limsup_{\e \to
0}\Norm{e^{i\frac{t}{2}\d_x^2}\W^\e}{L^6(\R^2)}\le 2\delta\, , 
\end{equation*}
we get
\begin{equation*}
\Norm{e^{i\frac{t}{2}\d_x^2}\phi}{L^6(\R^2)}\ge \frac{\beta}{2}\, , 
\end{equation*}
provided that $\delta>0$ is sufficiently small. This completes the
proof of Proposition~\ref{prop:2} in the general case. \\

When $(\U_0^\e)_\e$ satisfies \eqref{eq:echelle1}, there exists
$R=R(\delta)$ such that for every $\e$, 
\begin{equation*}
\Norm{\widehat{\U}_0^\e\1_{|\xi|\le R}}{L^2(\R)}\ge
\Norm{\widehat{\U}_0^\e}{L^2(\R)}  -\frac{\delta}{2}\, .
\end{equation*}
In the proof of Lemma~\ref{lem:1} (this is the step where the scales
$h^\e$ and cores in the Fourier side appear), we can therefore consider
$\widehat{\U}_0^\e\1_{|\xi|\le R}$ in place of
$\widehat{\U}_0^\e$. This implies that for any $j$,
$-\theta_j^\e/h_j^\e$ (the center of the balls we extract) and
$1/h_j^\e$ (the radius of the balls we extract) are uniformly
bounded. This means exactly that the sequence $(\xi_j^\e)_\e$ is
bounded for every $j$, and that $h_j^\e$ is bounded away from zero. As
mentioned already, up to an $\e$-independent dilation of the profiles
$\phi_j$, we deduce $ h^\e_j\ge 1$. 
\end{proof}

\begin{rema}\label{rema:dimension}
Why do we suppose $n=1$ or $2$ only? Essentially to have a refined
Strichartz estimate, as in \cite{MVV} in the case of space dimension
two, and in Proposition~\ref{prop:refined} for the one-dimensional
case. Notice that the proof uses the fact that
$2+\frac{4}{n}$ is an even integer, to decompose the
$L^{2+\frac{4}{n}}$ norm as a product. 
The restriction estimate
\eqref{eq:restriction} holds in higher dimensions. It is proved in 
\cite{Bourgain91} that if the space
dimension is $n\ge 3$, then such an 
estimate holds for some $q<2+\frac{4}{n}$, which is what
we use in the above computations (it holds more generally for
$q>2+\frac{4}{n-1}$, see \cite{TaoGAFA03}). 
\end{rema}

\section{Proof of Theorem~\ref{theo:profileNL}: nonlinear profile
decomposition} 
\label{sec:profileNL}

Roughly speaking, Theorem~\ref{theo:profileNL} is essentially a consequence of
Theorem~\ref{theo:profile} and of Strichartz 
inequalities, and is based on a perturbative analysis. This result has
no exact counterpart in \cite{MerleVega98}. Notice that one of the
key ingredients is Theorem~\ref{theo:profile}, and this is the only
reason why we have to restrict the space dimension. Since the approach
is very similar to 
\cite{BG3,Keraani01,IsabelleBSMF}, we shall only sketch the proof (see
\cite{KeraaniPhD} for more details). 

We prove the equivalence (i)$\Leftrightarrow$(ii); since the profiles
$\U_j$ are given by Theorem~\ref{theo:profile} and
Definition~\ref{def:profilNL}, and $r_\ell^\e$ is given by
Theorem~\ref{theo:profile}, only \eqref{eq:resteNL} has to be
proved. It follows from the perturbative argument of the proof
(i)$\Leftrightarrow$(ii). \\ 

\noindent \underline{(i)$\Rightarrow$(ii)}. Recall that $I_j^\e$ is
defined by $ I_j^\e := (h_j^\e)^{-2}\(I^\e-t_j^\e\)$, and that
$\U_j^\e$ is given by \eqref{eq:profNL}. We shall also denote
$\V_j^\e$ for the functions defined like in \eqref{eq:profNL}, with
$\U_j$ replaced by $e^{i\frac{t}{2}\Delta}\phi_j$, given by
Theorem~\ref{theo:profile}, that 
is $\V_j^\e= \H_j^\e (\phi_j)$ (see \eqref{eq:Hexpl}).

The function  $\rho_\ell^\e$ is  \emph{defined} by
$\rho_\ell^\e =\U^\e - \sum_{j=1}^\ell \U_j^\e - r_\ell^\e $. 
Denote $F(z)=\lambda |z|^{4/n}z$. The (expected) remainder $\rho_\ell^\e$
solves 
\begin{equation*}
i\d_t \rho_\ell^\e  +\frac{1}{2}\Delta \rho_\ell^\e  = f_\ell^\e\quad ;
\quad 
\rho_\ell^\e\big|_{t=0}= \sum_{j=1}^\ell\( \V_j^\e-\U_j^\e\)\big|_{t=0}\, ,
\end{equation*}
where 
\begin{equation*}
f_\ell^\e = F\Bigg( \rho_\ell^\e + \sum_{j=1}^\ell \U_j^\e + r_\ell^\e\Bigg) -
\sum_{j=1}^\ell  F\(\U_j^\e\)\, .
\end{equation*}
We use the orthogonality of the $\Gamma_j^\e$'s and the
assumption (i) to prove that \eqref{eq:resteNL} holds, that is:
$\dis\limsup_{\e \to 0} \(\Norm{\rho_\ell^\e}{L^\gamma(I^\e\times \R^n)}
+  \Norm{\rho_\ell^\e}{L^\infty(I^\e;L^2( \R^n))}\) \Tend \ell {+\infty}
0$. 

Once proved, this property implies (ii), since for some $\ell_0$
sufficiently large,
\begin{equation*}
\limsup_{\e \to 0 }\Norm{\U^\e}{L^\gamma(I^\e\times \R^n)}\leq
\sum_{j=1}^{\ell_0} \limsup_{\e \to 0
}\Norm{\U_j}{L^\gamma(I^\e_j\times \R^n)} +1<+\infty\, ,
\end{equation*}
by assumption (i). For $J^\e=[a^\e,b^\e]\subset I^\e$, Strichartz inequalities
yield 
\begin{equation*}
\Norm{\rho_\ell^\e}{L^\gamma(J^\e\times\R^n)}
+\Norm{\rho_\ell^\e}{L^\infty(J^\e;L^2(\R^n))}\lesssim
\Norm{\rho_\ell^\e(a^\e)}{L^2(\R^n)} + \Norm{f_\ell^\e}{L^{\gamma
'}(J^\e\times\R^n)}\, . 
\end{equation*}
From triangle and H\"older's inequalities,
\begin{align}
\Norm{f_\ell^\e}{L^{\gamma '}(J^\e\times\R^n)} \lesssim   
\Norm{\rho_\ell^\e}{L^{\gamma }}^\gamma &+ \Bigg\|\sum_{j=1}^\ell
\U_j^\e+r_\ell^\e\Bigg\|_{L^{\gamma }}^{\gamma -1}\Norm{\rho_\ell^\e}{L^{\gamma }}
\label{eq:321}\\
+ \Bigg\|\sum_{j=1}^\ell F(\U_j^\e) - F\Bigg(\sum_{j=1}^\ell \U_j^\e
\Bigg)\Bigg\|_{L^{\gamma  '}}& + \Bigg\|F\Bigg(\sum_{j=1}^\ell
\U_j^\e+r_\ell^\e\Bigg) - F\Bigg(\sum_{j=1}^\ell \U_j^\e 
\Bigg)\Bigg\|_{L^{\gamma  '}} \label{eq:123}.
\end{align}
The terms in \eqref{eq:123} are small by assumption (i), H\"older's
inequality and 
orthogonality (see for instance \cite{PG98}, and \cite{CFG} when
$\gamma$ is not an integer). The first term in
\eqref{eq:321} is treated by a bootstrap argument. We have to
take care of the second term in \eqref{eq:321}. 
The next lemma is proved in \cite{KeraaniPhD}. It allows to absorb this
linear term, thanks to a suitable partition of the interval $I^\e$. 
\begin{lem}\label{lem:absorb}
For every $\delta>0$, there exists an $\e$--dependent \emph{finite}
partition of $I^\e$,
\begin{equation*}
I^\e = \bigcup_{k=1}^{p(\delta)}J^\e_k\, ,
\end{equation*}
such that for every $1\le k\le p(\delta)$ and every $\ell \ge 1$,
\begin{equation*}
\limsup_{\e \to 0} \Bigg\|\sum_{j=1}^\ell \U_j^\e\Bigg\|_{L^\gamma (J^\e_k
\times\R^n)} \le \delta\, .
\end{equation*}
\end{lem}
\begin{proof}[Sketch of the proof]
By orthogonality, for every $\ell\ge 1$,
\begin{equation*}
\limsup_{\e \to 0} \Bigg\|\sum_{j=1}^\ell \U_j^\e\Bigg\|_{L^\gamma (I^\e
\times\R^n)}^\gamma = \sum_{j=1}^\ell \limsup_{\e \to 0}
\Norm{\U_j^\e}{L^\gamma (I^\e
\times\R^n)}^\gamma\, . 
\end{equation*}
On the other hand, the almost $L^2$--orthogonality \eqref{eq:orth}
and the conservation of mass for 
\eqref{eq:NLS} imply that for some $\ell(\delta)$,
\begin{equation*}
\Norm{\U_j}{L^2(\R^n)}\le \delta\, ,\quad \forall j\ge \ell(\delta)\, .
\end{equation*}
Using global existence results for small $L^2$ data (see
 e.g. \cite{CazCourant}), $\U_j$ is then defined globally in time,
 and from Strichartz estimate,
\begin{equation*}
\Norm{\U_j}{L^\gamma(\R\times\R^n)}\lesssim \Norm{\U_j}{L^2(\R^n)} =
 \Norm{\phi_j}{L^2(\R^n)} \, .
\end{equation*} 
Since $\gamma >2$ for any $n\ge 1$, we infer
\begin{equation*}
\sum_{j\ge \ell (\delta)}
\Norm{ \U_j}{L^\gamma (\R
\times\R^n)}^\gamma<+\infty\, . 
\end{equation*}
Using this and orthogonality, we infer
\begin{equation*}
\limsup_{\e \to 0}\Bigg\|\sum_{j=1}^\ell
\U_j^\e\Bigg\|_{L^\gamma (I^\e
\times\R^n)}^\gamma\le  \sum_{j=1}^{\ell(\delta)}\limsup_{\e \to
0} \Norm{\U_j^\e}{L^\gamma (I^\e \times\R^n)}^\gamma
+\frac{\delta}{2}\, . 
\end{equation*}
Thus, it suffices to construct a family of partial decompositions as
in the statement of the lemma, for every $1\le j\le \ell(\delta)$ and
such that 
\begin{equation*}
\limsup_{\e\to 0} \Norm{\U_j^\e}{L^\gamma(J_k^\e\times\R^n)}^\gamma
\le \frac{\delta}{2\ell(\delta)}\, ,\quad \forall 1\le k\leq
p(\delta)\, .
\end{equation*}
The final decomposition is obtained by intersecting all the partial
ones. 
Consider the case $j=1$, and denote by $I_1$ the maximal interval of
existence of $\U_1$. One checks that there exists a closed interval
$J_1$ such that 
\begin{equation*}
\limsup_{\e \to 0}I_1^\e = J_1\, ,\quad
\Norm{\U_1}{L^\gamma(J_1\times\R^n)}<+\infty\, .
\end{equation*}
We decompose $J_1$ as $ J_1 =\cup_{k=1}^{p_1(\delta)}J_{1k}$ so that
\begin{equation*}
\Norm{\U_1}{L^\gamma(J_{1k}\times\R^n)}<\frac{\delta}{2\ell(\delta)}\,
\quad \forall 1\le k \le p_1(\delta)\, .
\end{equation*}
At this first step, the intervals $J_1^\e$ are then obtained by scaling:
\begin{equation*}
J_k^\e = I^\e \cap  \( (h_1^\e)^{2} J_{1k} + t_1^\e\)\, .
\end{equation*}
Repeating this argument on each $J_k^\e $ a
finite number of times yields the lemma. 
\end{proof}
Choosing $\delta>0$ sufficiently small, Lemma~\ref{lem:absorb} allows
to prove that 
\begin{equation*}
\limsup_{\e \to 0} \(\Norm{\rho_\ell^\e}{L^\gamma(I^\e\times \R^n)}
+  \Norm{\rho_\ell^\e}{L^\infty(I^\e;L^2( \R^n))}\) \Tend \ell {+\infty}
0\, ,
\end{equation*}
thanks to an absorption argument for the linear term \eqref{eq:321},
orthogonality in 
the source term \eqref{eq:123}, and a bootstrap argument. \\

\noindent \underline{(ii)$\Rightarrow$(i)}. By assumption, there
exists $M>0$ such that 
\begin{equation*}
\limsup_{\e \to 0 }\Norm{\U^\e}{L^\gamma(I^\e\times
\R^n)} \le \frac{M}{2}\, .
\end{equation*}
Assume that (i) does not hold. Reorganizing the family of profiles, we
may assume that for some $\ell_0\ge 1$, $\U_j$ is not global -- that
is $\Norm{\U_j}{L^\gamma(\R\times\R^n)}=\infty$ -- if $1\le j\le
\ell_0$, and $\U_j$ is global for $j>\ell_0$. Indeed, if all the
profiles are defined globally in time, the problem is trivial. Thus,
we only have to consider a finite family of profiles, thanks to the
small data global existence results mentioned above.

Let $I_j$ denote the maximal interval of existence of $\U_j$, for
$1\le j\le \ell_0$. The failure of (i) means that there exists some
intervals $I_j(M)$ such that
\begin{equation*}
\frac{-t_j^\e}{(h_j^\e)^2} \in I_j(M)\subset I_j\cap I_j^\e \quad
\text{for }\e\ll 1\  
;\ M\le
\Norm{\U_j}{L^\gamma(I_j(M)\times\R^n)}<\infty\, ,\quad 1\le j\le
\ell_0\, .
\end{equation*}
Denote $I_j^\e(M)= (h_j^\e)^2 I_j(M) + t_j^\e$. Then $0 \in I_j^\e
(M)\subset I^\e$ for $\e \ll 1$ and 
\begin{equation}\label{eq:fim}
M\le
\limsup_{\e\to
0}\Norm{\U_j^\e}{L^\gamma(I_j^\e(M)\times\R^n)}<\infty\, ,\quad 1\le
j\le 
\ell_0\, .
\end{equation}
By permutation, extraction of a subsequence and considering the
backward and inward problems separately, we may take
\begin{equation*}
I_1^\e(M)\subset I_2^\e(M)\subset\ldots\subset I_{\ell_0}^\e(M)\, .
\end{equation*}
We infer
\begin{equation*}
\Norm{\U_j^\e}{L^\gamma(I_1^\e(M)\times\R^n)}<\infty\, ,\quad 1\le j\le
\ell_0\, .
\end{equation*}
We have 
\begin{equation}\label{eq:1159}
\limsup_{\e\to 0}\Norm{\U^\e}{L^\gamma(I_1^\e(M)\times\R^n)} \le
\limsup_{\e\to 0}\Norm{\U^\e}{L^\gamma(I^\e\times\R^n)}\le
\frac{M}{2}\, .
\end{equation}
Since $\U_j$ is global for $j>\ell_0$, (i) is satisfied with $I^\e$
replaced by $I_1^\e(M)$, and the first part of the proof yields
\eqref{eq:resteNL}. By orthogonality,
\begin{equation}\label{eq:1200}
\begin{aligned}
\limsup_{\e\to 0}\Norm{\U^\e}{L^\gamma(I_1^\e(M)\times\R^n)}^\gamma &=
\limsup_{\ell \to \infty}\Bigg( \limsup_{\e\to 0}\Big\|\sum_{j=1}^\ell
\U_j^\e \Big\|_{L^\gamma(I_1^\e(M)\times\R^n)}^\gamma\Bigg)\\
&=
\sum_{j=1}^\infty \limsup_{\e\to 0}\Norm{\U_j^\e
}{L^\gamma(I_1^\e(M)\times\R^n)}^\gamma\, . 
\end{aligned}
\end{equation}
In particular, \eqref{eq:1159} and \eqref{eq:1200} yield
\begin{equation*}
\limsup_{\e\to 0}\Norm{\U_1^\e}{L^\gamma(I_1^\e(M)\times\R^n)}\le
\frac{M}{2}\, , 
\end{equation*}
which contradicts \eqref{eq:fim}. Thus (i) holds, and we
saw in the first part of the proof that it implies \eqref{eq:resteNL}.

\section{Proof of Theorem \ref{theo:cns}: linearizability}
\label{sec:linear}
Using the scaling \eqref{eq:chgt}, we restate Theorem~\ref{theo:cns}. 
Define 
\begin{equation*}
\U_0^\e := \U^\e_{\mid t=0}\ \text{ and } V^\e :=
e^{i\frac{t}{2}\Delta}\U_0^\e \, . 
\end{equation*}
Then $u^\e$ and $\U^\e$ are
simultaneously linearizable on $I^\e$ in $L^2$. Moreover, $u^\e$ is
 linearizable on $I^\e$ in $H^1_\e$ if and only if   $\U^\e$ is
 linearizable on $I^\e$ in $H^1_{\sqrt\e}$. We now have to prove:
\begin{theo}\label{theo:cns'} Assume $n=1$ or $2$. Let $\U_0^\e$ 
bounded in $L^2(\R^n)$, $I^\e \ni 0$ a time interval. 
\begin{itemize}
\item $\U^\e$ is linearizable on $I^\e$ in $L^2$ if and only if
\begin{equation}\label{eq:linear'}
\limsup_{\e \to 0} \|\V^\e\|_{L^{2+4/n}(I^\e\times \R^n)}^{2+4/n} = 0
\, .
\end{equation}
\item Assume in addition that $\U_0^\e\in H^1$ and $\U_0^\e$ is
bounded in $H^1_{\sqrt\e}$. Then $\U^\e$ is linearizable on $I^\e$ in
$H^1_{\sqrt\e}$ if 
and only if \eqref{eq:linear'} holds. 
\end{itemize}
\end{theo}
\begin{proof}
We first prove that Condition \eqref{eq:linear'} is sufficient for
linearizability, thanks to the classical Strichartz estimates. In
particular, no restriction on the space dimension is necessary at this
stage. 
Denote $\W^\e = \U^\e -\V^\e$. It solves
\begin{equation}\label{eq:w}
     i\d_t \W^\e +\frac{1}{2}\Delta \W^\e
     =\lambda |\U^\e |^{4/n} \U^\e\, ,\quad
     \W^\e_{\mid t=0}=0\, .
\end{equation}
Since $\gamma=2+4/n$, we have 
\begin{equation*}
\frac{1}{\gamma'} = \frac{1}{\gamma} + \frac{4/n}{\gamma}\, .
\end{equation*}
Applying Strichartz estimate \eqref{eq:strichartzinhomogene} to
\eqref{eq:w}, along with H\"older's inequality, we have, for $t\in
I^\e$, 
\begin{equation*}
\begin{aligned}
\|\W^\e\|_{L^\gamma([0,t]\times \R^n)} &\lesssim 
\||\U^\e|^{4/n}\U^\e\|_{L^{\gamma '}([0,t]\times \R^n)} \lesssim 
\|\U^\e\|_{L^{\gamma }([0,t]\times \R^n)}^{1+\frac{4}{n}} \\
&\lesssim   \|\V^\e\|_{L^{\gamma }(I^\e\times \R^n)}^{1+\frac{4}{n}}+
\|\W^\e\|_{L^{\gamma }([0,t]\times \R^n)}^{1+\frac{4}{n}} \, .
\end{aligned}
\end{equation*}
Using Assumption \eqref{eq:linear'}, we apply a bootstrap argument: for
$\e$ sufficiently small, 
\begin{equation*}
\|\W^\e\|_{L^\gamma(I^\e\times \R^n)} \lesssim  \|\V^\e\|_{L^{\gamma
}(I^\e\times \R^n)}^{1+\frac{4}{n}}\, .
\end{equation*}
We infer that for $\e$ sufficiently small,
\begin{equation}\label{eq:upetit}
\|\U^\e\|_{L^\gamma(I^\e\times \R^n)} \lesssim  \|\V^\e\|_{L^{\gamma
}(I^\e\times \R^n)}+ \|\V^\e\|_{L^{\gamma
}(I^\e\times \R^n)}^{1+\frac{4}{n}}\, ,
\end{equation}
and \eqref{eq:linear'} holds with $\V^\e$ replaced by $\U^\e$. 
Applying the second part of Strichartz estimate
\eqref{eq:strichartzinhomogene} yields
\begin{equation*}
\|\W^\e\|_{L^\infty (I^\e;L^2(\R^n))} \lesssim
\||\U^\e|^{4/n}\U^\e\|_{L^{\gamma '}(I^\e\times \R^n)} \lesssim  
\|\U^\e\|_{L^{\gamma }(I^\e\times \R^n)}^{1+\frac{4}{n}}\Tend \e 0 0
\, , 
\end{equation*}
which is linearizability on $I^\e$ in $L^2$. 

Now assume that $\U^\e_0\in H^1(\R^n)$ is bounded in $H^1_{\sqrt\e}$. 
Differentiating
\eqref{eq:w} with respect to the space variable, we have
\begin{equation*}
\begin{aligned}
\|\sqrt\e \nabla_x \W^\e&\|_{L^\gamma(I^\e\times \R^n)} \lesssim 
\||\U^\e|^{4/n}\sqrt\e \nabla_x \U^\e\|_{L^{\gamma '}(I^\e\times
\R^n)} \\
&\lesssim 
\|\U^\e\|_{L^{\gamma }(I^\e\times \R^n)}^{4/n}\|\sqrt\e \nabla_x
\U^\e\|_{L^{\gamma }(I^\e\times \R^n)} \\ 
&\lesssim \|\U^\e\|_{L^{\gamma }(I^\e\times \R^n)}^{4/n} \left(
\|\sqrt\e \nabla_x \V^\e\|_{L^{\gamma }(I^\e\times \R^n)}+ 
\|\sqrt\e \nabla_x \W^\e\|_{L^{\gamma }(I^\e\times \R^n)} \right)\, .
\end{aligned}
\end{equation*}
From \eqref{eq:upetit} and \eqref{eq:linear'}, the term in $\sqrt\e
\nabla_x \W^\e$ on the 
right hand side can be absorbed by the left hand side for $\e$
sufficiently small. The uniform boundedness of $\sqrt\e \nabla_x
\V^\e$ in $L^\gamma(\R\times\R^n)$, which stems from the boundedness
of its data in L$^2$ and Strichartz estimate
\eqref{eq:strichartzhomogene}, shows that 
\begin{equation*}
\|\sqrt\e \nabla_x \W^\e\|_{L^\gamma(I^\e\times \R^n)}\Tend \e 0 0\, . 
\end{equation*}
Applying inhomogeneous Strichartz estimate
\eqref{eq:strichartzinhomogene} yields
\begin{equation*}
\|\sqrt\e \nabla_x \W^\e\|_{L^\infty(I^\e; L^2(\R^n))} \lesssim 
\|\U^\e\|_{L^{\gamma 
}(I^\e\times \R^n)}^{4/n} 
\|\sqrt\e \nabla_x \U^\e\|_{L^{\gamma }(I^\e\times \R^n)}\Tend
\e 0 0\, ,
\end{equation*}
which proves that $\U^\e$ is linearizable on $I^\e$ in
$H^1_{\sqrt\e}$. 
\medskip

We complete the proof of Theorem~\ref{theo:cns'} by showing that
Condition~\eqref{eq:linear'} is necessary for linearizability in
$L^2$ (hence for linearizability in $H^1_{\sqrt\e}$). The proof relies on
the profile decompositions stated in Theorems~\ref{theo:profile} and
\ref{theo:profileNL}. 
We consider two cases.\\

\noindent {\bf First case.} The family $(\U^\e)_{0<\e\leq 1}$ is
uniformly bounded in $L^\gamma(I^\e\times\R^n)$. \\
In that case, we can use Theorems~\ref{theo:profile} and
\ref{theo:profileNL} to deduce the following lemma. 
\begin{lem}\label{lem:decompdifference}
Assume $n=1$ or $2$. Let $\U^\e_0$ bounded in $L^2(\R^n)$, $I^\e =[0,T^\e[
$ a (possibly unbounded) time interval, and assume
that $\U^\e$ is bounded in $L^\gamma(I^\e\times\R^n)$. Then up to the
extraction of a subsequence, there exist an 
orthogonal family $(h^\e_j,t^\e_j,x^\e_j,\xi^\e_j)_{j \in \N}$ 
in $ \R_+\setminus\{0\}\times\R\times \R^n\times \R^n$ and a family
$(\phi_j)_{j \in \N}$ bounded in~$L^2(\R^n) $, such that if $\V_j =
e^{i\frac{t}{2}\Delta}\phi_j$ and $\U_j$ is given by
Definition~\ref{def:profilNL}, we have:
\begin{equation}\label{eq:serie}
\limsup_{\e \to 0}\Norm{ \U^\e
-\V^\e}{L^\gamma(I^\e\times\R^n)}^\gamma = \sum_{j=1}^\infty
\limsup_{\e \to 0}\left\| \U_j
-\V_j\right\|_{L^\gamma(I^\e_j\times\R^n)}^\gamma\, ,
\end{equation}
where $I_j^\e = (h_j^\e)^{-2}(I^\e -t_j^\e)$. 
In addition, for every fixed $\e>0$, none of the terms in the series
is zero. 
\end{lem}
\begin{proof}[Proof of Lemma~\ref{lem:decompdifference}]
From Theorems~\ref{theo:profile} and \ref{theo:profileNL}, there exist
an orthogonal 
family $(h^\e_j,t^\e_j,x^\e_j,\xi^\e_j)_{j \in \N}$  
in $ \R_+\setminus\{0\}\times\R\times \R^n\times \R^n$ and a family
$(\phi_j)_{j \in \N}$ bounded in~$L^2(\R^n) $, such that if $\V_j =
e^{i\frac{t}{2}\Delta}\phi_j$ and $\U_j$ is given by
Definition~\ref{def:profilNL} (up to the extraction of a subsequence),
we have, for any $\ell \in \N$, 
\begin{equation}\label{eq:decompdifference}
\U^\e(t,x)
-\V^\e(t,x) = \sum_{j=1}^\ell \H_j^\e\left(
 \U_j\big|_{t=0}
-\V_j\big|_{t=0}\right)(t,x) + \rho^\e_\ell(t,x)\, ,
\end{equation}
with $\limsup_{\e\to 0}\|\rho^\e_\ell\|_{L^\gamma(I^\e\times \R^n)}\to 0$
as $\ell \to +\infty$. The scales $h_j^\e$, cores
$(t^\e_j,x^\e_j,\xi^\e_j)$ and initial profiles $\phi_j$ are the same
for $\U^\e$ and $\V^\e$, since
they are given by the profile decomposition for the initial data
$\U^\e_{\mid t=0} = \V^\e_{\mid t=0} =\U^\e_{0}$.  
Since the family
$(h^\e_j,t^\e_j,x^\e_j,\xi^\e_j)_{j \in \N}$ is orthogonal, we have,
for any $\ell$, 
\begin{align*}
\limsup_{\e \to 0}\left\| \U^\e
-\V^\e\right\|_{L^\gamma(I^\e\times\R^n)}^\gamma = &\sum_{j=1}^\ell
\limsup_{\e \to 0}\left\| \U_j
-\V_j\right\|_{L^\gamma(I_j^\e\times\R^n)}^\gamma\\
&+\limsup_{\e\to
 0}\|\rho^\e_\ell\|_{L^\gamma(I^\e\times \R^n)}^\gamma   .
\end{align*}
Letting $\ell \to +\infty$ yields \eqref{eq:serie}. Now assume that
for a fixed $\e>0$, one of the terms in the series \eqref{eq:serie} is
zero. This means that two solutions of the nonlinear Schr\"odinger
equation \eqref{eq:NLS} and of the Schr\"odinger equation respectively
coincide on
the non-trivial time interval $I_{j_0}^\e$. Uniqueness for these two
equations shows 
that necessarily $\U_{j_0}=\V_{j_0}\equiv 0$, in which case the family
$(\U_j,\V_j)_j$ can be relabeled to avoid null terms. 
\end{proof}

\begin{defin}
Let $\delta^\e>0$ and $a^\e\in \R$. We say that the interval
$]a^\e,a^\e+\delta^\e[$ is 
\emph{asymptotically trivial} in either of the following cases:
\begin{itemize}
\item $a^\e \to +\infty$ as $\e \to 0$, or
\item $a^\e +\delta^\e\to -\infty$ as $\e \to 0$, or
\item $\delta^\e\to 0$ as $\e \to 0$. 
\end{itemize}
\end{defin}

\begin{lem}\label{lem:petit}
Under the assumptions of Lemma~\ref{lem:decompdifference}, if $\|\U^\e
-\V^\e\|_{L^\gamma (I^\e\times\R^n)}\to 0$ as $\e \to 0$, then
$\|\V^\e\|_{L^\gamma (I^\e\times\R^n)}\to 0$ as $\e \to 0$. 
\end{lem}
\begin{proof}[Proof of Lemma~\ref{lem:petit}]
From Lemma~\ref{lem:decompdifference}, if $\|\U^\e
-\V^\e\|_{L^\gamma (I^\e\times\R^n)}\to 0$, then every interval
$I_j^\e$ is asymptotically trivial. The
profile decomposition for $\V^\e$ yields
\begin{equation*}
\V^\e(t,x) = \sum_{j=1}^\ell \H_j^\e\left(
\V_j\right)(t,x) + r^\e_\ell(t,x)\, ,
\end{equation*}
with $\limsup_{\e\to 0}\|r^\e_\ell\|_{L^\gamma(I^\e\times \R^n)}\to 0$
as $\ell \to +\infty$. Fix $\ell \in \N$. We infer from the
orthogonality of $(h^\e_j,t^\e_j,x^\e_j,\xi^\e_j)_{j \in \N}$ that
\begin{equation*}
\limsup_{\e\to 0}\|\V^\e\|_{L^\gamma (I^\e\times\R^n)}^\gamma =
\sum_{j=1}^\ell \limsup_{\e \to 0}\|
\V_j \|_{L^\gamma (I_j^\e\times\R^n)}^\gamma
+
\limsup_{\e \to 0}\|r^\e_\ell \|_{L^\gamma (I^\e\times\R^n)}^\gamma\, .
\end{equation*}
Since all the intervals $I_j^\e$ are
asymptotically trivial, every term in the sum is zero, and we have
\begin{equation*}
\limsup_{\e\to 0}\|\V^\e\|_{L^\gamma (I^\e\times\R^n)} =
\limsup_{\e \to 0}\|r^\e_\ell \|_{L^\gamma (I^\e\times\R^n)}\, .
\end{equation*}
Since the left hand side is independent of $\ell$, we conclude that
both terms are zero, which completes the proof of
Lemma~\ref{lem:petit}. 
\end{proof}

We can now complete the proof of Theorem~\ref{theo:cns'} in the case
where the family $(\U^\e)_{0<\e\leq 1}$ is
uniformly bounded in $L^\gamma(I^\e\times\R^n)$. Assume that $\U^\e$
is linearizable on $I^\e$ in $L^2$. From Lemma~\ref{lem:petit}, it is
enough to prove that 
\begin{equation*}
 \|\U^\e-\V^\e\|_{L^\gamma (I^\e\times\R^n)}\Tend \e 0 0\, .
\end{equation*}
If it were not so, then from Lemma~\ref{lem:decompdifference}, there
would exist $j_0$ such that the interval $I_{j_0}^\e$ is not
asymptotically trivial. Up to the 
extraction of a subsequence, we can assume that there exist $a<b$
independent of $\e$ and $\e_0$ such that for $0<\e\leq \e_0$, 
\begin{equation*}
]a,b[\subset I_{j_0}^\e=\left[-\frac{t_{j_0}^\e}{(h_{j_0}^\e)^2}
,\frac{T^\e-t_{j_0}^\e}{(h_{j_0}^\e)^2}\right[\, .
\end{equation*}
Let $\ell > j_0$. Apply the operator $(\H_{j_0}^\e)^{-1}$ to
\eqref{eq:decompdifference}, and take the weak limit in ${\mathcal
D}'(]a,b[\times \R^n)$. By orthogonality,
\begin{equation}\label{eq:limitefaible}
w\!\!-\!\!\lim (\H_{j_0}^\e)^{-1}(\U^\e -\V^\e) = (\U_{j_0}-\V_{j_0}){\bf
1}_{]a,b[}(t) + w\!\!-\!\!\lim (\H_{j_0}^\e)^{-1}\rho^\e_\ell\, .
\end{equation}
Denote $\w_\ell := w\!\!-\!\!\lim (H_{j_0}^\e)^{-1}\rho^\e_\ell$. We
have
\begin{equation*}
\|\w_\ell\|_{L^\gamma (]a,b[\times \R^n)}\leq \liminf_{\e \to 0}
\|(\H_{j_0}^\e)^{-1}\rho^\e_\ell \|_{L^\gamma (]a,b[\times \R^n)} \leq
\liminf_{\e \to 0} \|\rho^\e_\ell \|_{L^\gamma (I^\e\times \R^n)}\Tend
\ell {+\infty} 0\, .
\end{equation*}
In \eqref{eq:limitefaible}, $\w_\ell$ is the only term possibly
depending on $\ell$, therefore it is zero, and 
$w\!\!-\!\!\lim (\H_{j_0}^\e)^{-1}(\U^\e -\V^\e) \not = 0$. Since
$\H_{j_0}^\e$ is unitary on $L^2(\R^n)$, we have, for $0<\e\leq \e_0$, 
\begin{equation*}
\left\| (\H_{j_0}^\e)^{-1}(\U^\e
-\V^\e)\right\|_{L^\infty(]a,b[;L^2(\R^n))} \leq \left\|
\U^\e 
-\V^\e\right\|_{L^\infty(I^\e;L^2(\R^n))}\, .
\end{equation*}
The right hand side goes to zero as $\e \to 0$ since $\U^\e$ is
assumed to be linearizable on $I^\e$ in $L^2$. Therefore the left hand
side goes to zero. This is impossible, since the weak limit is not
zero. This contradiction shows that we can apply
Lemma~\ref{lem:petit}, and complete the proof of
Theorem~\ref{theo:cns'} in the case 
where the family $(\U^\e)_{0<\e\leq 1}$ is
uniformly bounded in $L^\gamma(I^\e\times\R^n)$.

\noindent {\bf Second case.} There exists a subsequence of
$(\U^\e)_{0<\e\leq 1}$, still denoted $\U^\e$, such that 
\begin{equation*}
\left\| \U^\e\right\|_{L^\gamma(I^\e\times\R^n)}\Tend \e 0 +\infty\, .
\end{equation*}
Then there exists $\tau^\e \in I^\e$ such that for every $\e\in ]0,1]$, 
\begin{equation}\label{eq:satur}
\left\| \U^\e\right\|_{L^\gamma([0,\tau^\e[\times\R^n)}=1\, .
\end{equation}
We can mimic the proof of the first case on the time interval
$[0,\tau^\e[$. Lemma~\ref{lem:petit} shows that $\left\|
\U^\e\right\|_{L^\gamma([0,\tau^\e[\times\R^n)}\to 0$ as $\e \to 0$,
which contradicts \eqref{eq:satur}. Therefore the second case never
occurs, and the proof of Theorem~\ref{theo:cns'} is complete.  
\end{proof}

\begin{rema*} The above proof of linearizability relies on the
profile decompositions (linear and nonlinear). Note that in
\cite{CFG}, the proof of 
linearizability used only the conservations of mass and energy, and
Strichartz inequalities. Only after the linearizability
criterion had been proved, a (linear) profile decomposition was used. 
\end{rema*}

\section{Obstructions to linearizability}
\label{sec:decomp}

\subsection{Profile decomposition}

In this paragraph, we show how to deduce Corollary~\ref{cor:nonlin}
from Theorem~\ref{theo:profile}. 

Resuming the scaling \eqref{eq:chgt}, \eqref{eq:DACI} is exactly the
result given by the first part of Theorem~\ref{theo:profile} on the
time interval $[0,T]$ when considering the trace
$t=0$. We use the first part of Theorem~\ref{theo:profile} because 
Theorem~\ref{theo:cns} reduces our problem to the study of a
solution to the \emph{linear} Schr\"odinger equation. 
Notice that even if
we considered a defocusing  nonlinearity
($\lambda =+1$), with $u_0^\e$ bounded in
$H^1_\e$, we could not claim that $\U^\e$ is uniformly bounded in
$L^\gamma([0,T]\times \R^n)$.
This is because we do not know that $H^1$ solutions to
\eqref{eq:NLS} with $\lambda =+1$ decay like solutions to the free
equations as time goes to infinity (this is known in $\Sigma$); this
issue is related to the asymptotic 
completeness of wave operators in $H^1$. 

Working with the functions $\U^\e$ and $\V^\e$, \eqref{eq:DACI} writes:
\begin{equation}\label{eq:DACI'}
\begin{aligned}
&\U^\e_0(x) =\sum_{j=1}^\ell \widetilde \H_j^\e (\phi_j)(x)
+\W^\e_\ell(x)\, ,&\\
\text{with }&
\limsup_{\e \to 0}\left\|e^{i\frac{t}{2}\Delta}
\W^\e_\ell\right\|_{L^{2+4/n}(\R\times\R^n)}^{2+4/n}
\Tend \ell {+\infty} 0\, . &
\end{aligned}
\end{equation}
From \eqref{eq:Hexpl},
\begin{equation*}
\left\|e^{i\frac{t}{2}\Delta} \widetilde \H_j^\e
(\phi_j)\right\|_{L^\gamma ([0,T]\times \R^n)} = \|
\V_j\|_{L^\gamma (I^\e_j \times \R^n)}\, , \text{ with }I^\e_j = \left[
\frac{-t^\e_j}{(h^\e_j)^{2}} \virgp
 \frac{T-t^\e_j}{(h^\e_j)^{2}}\right]\, .
\end{equation*}
If $I^\e_j$ is asymptotically trivial, then $ \widetilde \H_j^\e
(\phi_j)$ can be incorporated into the
remainder term $\W^\e_\ell$, a case which can be excluded, up to
relabeling our family of sequences. This means that we can assume:
\begin{equation*}
\frac{-t^\e_j}{(h^\e_j)^{2}} \not\to +\infty\ ,\quad
\frac{T-t^\e_j}{(h^\e_j)^{2}} \not \to -\infty\ ,\quad\text{and }\ 
\frac{T}{(h^\e_j)^{2}}\not \to 0\, .
\end{equation*}
The first two points imply the properties on $t^\e_j$ stated in
Corollary~\ref{cor:nonlin}. We infer from the last point that $h^\e_j$
is bounded, by $1$ up to the extraction of a
subsequence and an $\e$-independent dilation of the profiles
$\phi_j$.\\  

Now suppose that $u_0^\e\in H^1$ and is bounded in $H^1_\e$. 
Then for every $j$, $\xi^\e_j = \O(\e^{-1/2})$ as
$\e\to 0$. To see this, introduce the scaling
\begin{equation}\label{eq:scal}
\psi^\e(t,x)= \e^{n/2}u^\e(\e t,\e x)\, .
\end{equation}
The function $\psi^\e$ solves \eqref{eq:NLS}, and the family
$(\psi^\e(0,x))_{0<\e \leq 1}$ is bounded in $H^1$. The plane
oscillations in the decomposition for $\psi^\e$ are, from
\eqref{eq:scal}, $e^{ix\cdot \xi^\e_j \sqrt\e}$. From the second point
of Theorem~\ref{theo:profile}, we infer that
$\xi^\e_j\sqrt\e =\O(1)$. We also deduce the lower bound $h^\e_j\geq
\sqrt \e$.


Finally, \eqref{eq:DA} is obtained from \eqref{eq:DACI} \emph{via} the
classical formula (see e.g. \cite{Rauch91})
\begin{equation}\label{eq:DASchrodinfini}
e^{i\frac{t}{2}\Delta}\phi =
e^{in\frac{\pi}{4}}e^{i\frac{x^2}{2t}}\frac{1}{|t|^{n/2}} \widehat\phi
\left(\frac{x}{t}\right)+o(1)\quad \text{in }L^2(\R^n),\text{ as }t\to
-\infty\, .
\end{equation}

\subsection{Nonlinear superposition}

We know assume $\lambda =+1$. 
The decomposition \eqref{eq:DACI} is necessary for the nonlinear term
in \eqref{eq:ck} to have a leading order influence on finite term
intervals. The aim of this section is to provide an argument
suggesting that it is sufficient. As mentioned before, the gap between belief
and proof is related to the asymptotic completeness of wave operators
in $H^1$. 

Suppose the initial data $u^\e_0$ has the form \eqref{eq:DACI} for a
fixed $\ell$ and a linearizable remainder: there exists $T>0$ such that
\begin{equation}\label{eq:CI}
\begin{aligned}
&u^\e_0(x) =\sum_{j=1}^\ell \widetilde H_j^\e (\phi_j)(x)
+w^\e(x)\, ,&\\
\text{where }&
\widetilde H_j^\e (\phi_j)(x) =e^{ix\cdot
 \xi^\e_j/\sqrt\e}e^{-i\e\frac{t^\e_j}{2}\Delta}\left(
\frac{1}{(h^\e_j\sqrt\e)^{n/2}}\phi_j
 \left(\frac{x-x^\e_j}{h^\e_j\sqrt\e} \right)\right)\, ,&\\
\text{and }&
\limsup_{\e \to 0}\e\|e^{i\e\frac{t}{2}\Delta}
w^\e\|_{L^{2+4/n}([0,T]\times\R^n)}^{2+4/n}= 0\, . &
\end{aligned}
\end{equation}
If we assume that $\phi_j\in\Sigma$ for every $j\in
\{1,\ldots,\ell\}$, then we can take advantage of the global
well-posedness and the existence of a complete scattering theory 
for \eqref{eq:NLS} in $\Sigma$ when $\lambda=+1$. Moreover, we may
assume that $t^\e_j/(h^\e_j)^2$ converges as $\e\to 0$ for every $j$. 
Let $v^\e_j$ be the
solution of the initial value problem
\begin{equation}\label{eq:vj}
i\e \d_t v^\e_j +\frac{1}{2}\e^2\Delta v^\e_j = \e^2
|v^\e_j|^{4/n}v^\e_j\quad
 ;\quad v^\e_{j\mid t=0} = \widetilde H_j^\e (\phi_j)\,.
\end{equation}
For every $j$, the following asymptotics holds in $L^\infty(\R;L^2)$
as $\e$ goes to zero:
\begin{equation*}
v^\e_j(t,x) = e^{i\frac{x\cdot
 \xi^\e_j}{\sqrt\e}-i\frac{t}{2}(\xi^\e_j)^2}
\frac{1}{(h^\e_j\sqrt\e)^{n/2}}\V_j
 \left(\frac{t-t^\e_j}{(h^\e_j)^2}\virgp
\frac{x-x^\e_j-t\xi^\e_j}{h^\e_j\sqrt\e}\right)\, +o(1),
\end{equation*}
where $\V_j$ is given by
\begin{equation*}
i\d_t \V_j +\frac{1}{2}\Delta \V_j = |\V_j|^{4/n}\V_j\quad ;\quad
e^{-i\frac{t}{2}\Delta}\V_j(t)\big|_{t=-\lim
t^\e_j/(h^\e_j)^2}=\phi_j\, .
\end{equation*}
Notice that the above problem may be an initial value problem or a
scattering problem according to the value of $\lim
t^\e_j/(h^\e_j)^2$. We see that $v^\e_j$ has a genuine nonlinear
behavior on $[0,T]$ if $\lim
t^\e_j/(h^\e_j)^2\not = -\infty$ and $\lim
(T-t^\e_j)/(h^\e_j)^2\not = -\infty$ (compare with
Corollary~\ref{cor:nonlin}). 

Following the lines of \cite{BG3} and \cite{CFG}, the next result
can be shown, thanks to the linearizability criterion given by
Theorem~\ref{theo:cns}. We leave out the proof here, for it bears no
new idea. 
\begin{theo}\label{theo:superp}
Assume $n=1$ or $2$, $\lambda=+1$, and let $u_0^\e$ be given by
\eqref{eq:CI} with $\phi_j \in \Sigma$ and an 
orthogonal family $(h^\e_j,t^\e_j,x^\e_j,\xi^\e_j)_{j \in \N}$ such
that $\sqrt\e \leq h^\e_j\leq 1$.
Then the following asymptotics holds in
$L^\infty([0,T];L^2)$ as $\e$ goes to zero,
\begin{equation*}
u^\eps = \sum_{j=1}^\ell v_j^\e + e^{i\e\frac{t}{2}\Delta}
w^\e +o(1)\, , 
\end{equation*}
where each $v^\e_j$ solves \eqref{eq:vj}. 
\end{theo}

\section{Blowing up solutions}
\label{sec:blowup}
Assume $n=1$ or $2$. Let $\U$ be an
$L^2$-solution to \eqref{eq:NLS} which blows up\footnote{The general
consensus is that even in the $L^2$ framework, this can occur
only in the attractive case $\lambda <0$}  at time $T>0$ (not
before),
\begin{equation}\label{eq:defblowup}
\int_0^T\!\!\!\! \int_{\R^n}|\U(t,x)|^{2+\frac{4}{n}}dxdt =+\infty\, .
\end{equation}
Let $(t_k)_{k\in \N}$ be a sequence going to $T$ as $k \to
+\infty$, with $t_k<T$ for every $k$. Denote $\e_k =T-t_k$, and define
\begin{equation*}
u^\e(t,x)= \U(\e t+T-\e,x)\, ,
\end{equation*}
where the notation $\e$ stands for $\e_k$. Then $u^\e$ solves
\eqref{eq:ck}. The function $\U$ blows up at time $T$  if and only if
$u^\e$ is not linearizable on 
$[0,1]$ (in $L^2$), from Theorem~\ref{theo:cns} and its proof. 
The function $v^\e$ is given by 
\begin{equation*}
v^\e(t,x)=e^{i\e\frac{t}{2}\Delta}u_0^\e(x) =
e^{i\e\frac{t}{2}\Delta}\U(T-\e,x)\, .
\end{equation*}
Define 
\begin{equation*}
\V^\e(t,x)= v^\e\left( \frac{t-T}{\e}+1,x\right)\, .
\end{equation*}
Since $u^\e$ is not linearizable on $[0,1]$, we have
$\dis \liminf_{\e \to 0}\e\|v^\e\|_{L^\gamma ([0,1]\times\R^n)}^\gamma >0$.
From Corollary~\ref{cor:nonlin},  up to the 
extraction of a subsequence, 
\begin{equation*}
\begin{aligned}
&u^\e_0(x) =\sum_{j=1}^\ell \widetilde H_j^\e (\phi_j)(x)
+w^\e_\ell(x)\, ,&\\
\text{where }&
\widetilde H_j^\e (\phi_j)(x) =e^{ix\cdot
 \xi^\e_j/\sqrt\e}e^{-i\e\frac{t^\e_j}{2}\Delta}\left(
\frac{1}{(h^\e_j\sqrt \e)^{n/2}}\phi_j
 \left(\frac{x-x^\e_j}{h^\e_j\sqrt\e} \right)\right)\, ,&\\
\text{and }&
\limsup_{\e \to 0}\e\|e^{i\e\frac{t}{2}\Delta}
w^\e_\ell\|_{L^{2+4/n}(\R\times\R^n)}^{2+4/n}
\Tend \ell {+\infty} 0\, . &
\end{aligned}
\end{equation*}
Recall that from \eqref{eq:Hexpl},
\begin{equation*}
\widetilde H_j^\e (\phi_j)(x) =e^{ix\cdot \xi_j^\e/\sqrt\e}
  \frac{1}{(h^\e_j\sqrt\e)^{n/2}} \V_j \left(
\frac{-t^\e_j}{(h^\e_j)^{2}}\virgp \frac{x-x^\e_j}{h^\e_j\sqrt\e}\right)
, \
\text{where }\V_j(t)=e^{i\frac{t}{2}\Delta} \phi_j\, .
\end{equation*}
Moreover, we can assume 
\begin{equation*}
\frac{-t^\e_j}{(h^\e_j)^{2}} \not\to +\infty\ ,\quad
\frac{1-t^\e_j}{(h^\e_j)^{2}} \not \to -\infty\ ,\quad\text{and }\ 
\frac{1}{(h^\e_j)^{2}}\not \to 0\, ,
\end{equation*}
for otherwise, the corresponding profile may be incorporated into the
remainder $w^\e_\ell$. This implies that for every $j$,
$(h^\e_j)_{j\in \N}$ and $(t^\e_j)_{j\in \N}$ are  bounded
sequences. Up to extracting a subsequence, we distinguish two cases:
\begin{equation*}
\frac{t^\e_j}{(h^\e_j)^2}\to \lambda \in \R\text{ as }\e\to 0\, ,\quad
\text{or}\quad
\frac{t^\e_j}{(h^\e_j)^2}\to +\infty\, .
\end{equation*}
In the first case, we set  
$y^k_j=\xi^\e_j/\sqrt\e$, $x^k_j=x^\e_j$,  
$\rho^k_j = 
h^\e_j\sqrt\e=h^\e_j\sqrt{T-t_k}\leq \sqrt{T-t_k}$ and $\U_j =
\V_j(-\lambda)$. 
In the second case, we infer from
\eqref{eq:DASchrodinfini} that in $L^2$, 
\begin{equation*}
\widetilde H_j^\e (\phi_j)(x) \Eq \e 0 e^{in\frac{\pi}{4}+ix\cdot
\xi^\e_j/\sqrt\e} e^{-i\frac{|x-x^\e_j|^2}{2\e t^\e_j}}\left(
\frac{h^\e_j}{t^\e_j\sqrt\e} \right)^{n/2} \widehat \phi_j \left( 
\frac{h^\e_j}{t^\e_j\sqrt\e}(x-x^\e_j) \right)\, .
\end{equation*}
We set $y^k_j=\xi^\e_j/\sqrt\e$, $x^k_j=x^\e_j$,
and $\widetilde\U_j =e^{in\pi/4}\widehat \phi_j$, and the proof of
Corollary~\ref{cor:blowup} is complete, up to relabeling the family of
sequences and possibly taking some $\U_j$ or some $\widetilde\U_j$
equal to zero. 

\begin{rema*}
When  only one profile is present, quadratic oscillations are not
relevant near the blow-up time. Assume
\begin{equation*}
\begin{aligned}
&u^\e_0(x) =\widetilde H^\e (\phi)(x)
+w^\e(x)\, ,&\\
\text{where }&
\widetilde H^\e (\phi)(x) =e^{ix\cdot
 \xi^\e/\sqrt\e}e^{-i\e\frac{t^\e}{2}\Delta}\left(
\frac{1}{(h^\e\sqrt\e)^{n/2}}\phi
 \left(\frac{x-x^\e}{h^\e\sqrt\e} \right)\right)\, ,&\\
\text{and }&
\limsup_{\e \to 0}\e\|e^{i\e\frac{t}{2}\Delta}
w^\e\|_{L^{2+4/n}(\R\times\R^n)}^{2+4/n}
\Tend \ell {+\infty} 0\, . &
\end{aligned}
\end{equation*}
Since there is blow-up at time $T$,
\begin{equation}\label{eq:nlapres}
\liminf_{\e \to 0}\e\|v^\e\|_{L^\gamma ([0,1]\times\R^n)}^\gamma >0\, . 
\end{equation}
On the other hand, we also have
\begin{equation}\label{eq:nlavant}
\liminf_{\e \to 0}\|\V^\e\|_{L^\gamma ([0,T-\e]\times\R^n)}^\gamma = 
\liminf_{\e \to 0}\e\|v^\e\|_{L^\gamma ([1-T/\e,0]\times\R^n)}^\gamma
>0\, .
\end{equation}
If this limit
was zero, then $v^\e$ would be linearizable in $L^2$ on $[1-T/\e,0]$,
and
\begin{equation*}
\liminf_{\e \to 0}\e\|u^\e\|_{L^\gamma ([1-T/\e,0]\times\R^n)}^\gamma
=0=\liminf_{\e \to 0}\|\U\|_{L^\gamma ([0,T-\e]\times\R^n)}^\gamma\, 
 , 
\end{equation*}
which contradicts \eqref{eq:defblowup}. Recall
\begin{equation*}
\widetilde H^\e (\phi)(x) =e^{ix\cdot \xi^\e
/\sqrt\e}  \frac{1}{(h^\e\sqrt\e)^{n/2}} \V \left(
\frac{-t^\e}{(h^\e_j)^{2}}\virgp \frac{x-x^\e}{h^\e\sqrt\e}\right)
, \
\text{where }\V(t)=e^{i\frac{t}{2}\Delta} \phi\, .
\end{equation*}
From \eqref{eq:nlapres}, we have
$ -t^\e/(h^\e_j)^{2} \not \to +\infty$, 
and from \eqref{eq:nlavant},
$-t^\e/(h^\e_j)^{2} \not \to -\infty$.
Therefore, up to an extraction, $-t^\e/(h^\e_j)^{2} \to
\lambda \in \R$, and we are left with a profile only, and no quadratic
oscillation. 
\end{rema*}


\noindent {\bf Acknowledgments}. The authors are grateful to Clotilde
Fermanian--Kammerer and Isabelle Gallagher for stimulating discussions on
this work.

\bibliographystyle{amsplain}
\bibliography{carles}

\providecommand{\bysame}{\leavevmode\hbox to3em{\hrulefill}\thinspace}
\providecommand{\MR}{\relax\ifhmode\unskip\space\fi MR }
\providecommand{\MRhref}[2]{%
  \href{http://www.ams.org/mathscinet-getitem?mr=#1}{#2}
}
\providecommand{\href}[2]{#2}
\begin{thebibliography}{10}

\bibitem{BG3}
H.~Bahouri and P.~G{\'e}rard, \emph{High frequency approximation of solutions
  to critical nonlinear wave equations}, Amer. J. Math. \textbf{121} (1999),
  no.~1, 131--175. \MR{2000i:35123}

\bibitem{Bourgain91}
J.~Bourgain, \emph{Besicovitch type maximal operators and applications to
  {F}ourier analysis}, Geom. Funct. Anal. \textbf{1} (1991), no.~2, 147--187.

\bibitem{Bourgain95}
\bysame, \emph{Some new estimates on oscillatory integrals}, Essays on Fourier
  analysis in honor of Elias M. Stein (Princeton, NJ, 1991), Princeton Math.
  Ser., vol.~42, Princeton Univ. Press, Princeton, NJ, 1995, pp.~83--112.
  \MR{96c:42028}

\bibitem{Bourgain98}
\bysame, \emph{Refinements of {S}trichartz' inequality and applications to
  {$2$}{D}-{NLS} with critical nonlinearity}, Internat. Math. Res. Notices
  (1998), no.~5, 253--283.

\bibitem{BourgainWang}
J.~Bourgain and W.~Wang, \emph{Construction of blowup solutions for the
  nonlinear {S}chr\"odinger equation with critical nonlinearity}, Ann. Scuola
  Norm. Sup. Pisa Cl. Sci. (4) \textbf{25} (1997), no.~1--2, 197--215.

\bibitem{Ca2}
R.~Carles, \emph{Geometric optics with caustic crossing for some nonlinear
  {S}chr\"odinger equations}, Indiana Univ. Math. J. \textbf{49} (2000), no.~2,
  475--551.

\bibitem{CFG}
R.~Carles, C.~Fermanian, and I.~Gallagher, \emph{On the role of quadratic
  oscillations in nonlinear {S}chr{\"o}dinger equations}, J. Funct. Anal.
  \textbf{203} (2003), no.~2, 453--493.

\bibitem{CazCourant}
T.~Cazenave, \emph{Semilinear {S}chr\"odinger equations}, Courant Lecture Notes
  in Mathematics, vol.~10, New York University Courant Institute of
  Mathematical Sciences, New York, 2003.

\bibitem{CW89}
T.~Cazenave and F.~Weissler, \emph{Some remarks on the nonlinear
  {S}chr\"odinger equation in the critical case}, Lect. Notes in Math., vol.
  1394, Springer-Verlag, Berlin, 1989, pp.~18--29.

\bibitem{CW90}
\bysame, \emph{The {C}auchy problem for the nonlinear {S}chr\"odinger equation
  in ${H}^s$}, Nonlinear Anal. TMA \textbf{14} (1990), 807--836.

\bibitem{Fefferman83}
C.~Fefferman, \emph{The uncertainty principle}, Bull. Amer. Math. Soc. (N.S.)
  \textbf{9} (1983), no.~2, 129--206.

\bibitem{IsabelleBSMF}
I.~Gallagher, \emph{Profile decomposition for solutions of the
  {N}avier-{S}tokes equations}, Bull. Soc. Math. France \textbf{129} (2001),
  no.~2, 285--316. \MR{2002h:35235}

\bibitem{PG98}
P.~G{\'e}rard, \emph{Description du d\'efaut de compacit\'e de l'injection de
  {S}obolev}, ESAIM Control Optim. Calc. Var. \textbf{3} (1998), 213--233
  (electronic). \MR{99h:46051}

\bibitem{GV82}
J.~Ginibre and G.~Velo, \emph{Sur une \'equation de {S}chr\"odinger non
  lin\'eaire avec interaction non locale}, Nonlinear partial differential
  equations and their applications, {C}oll\`ege de {F}rance {S}eminar
  (H.~Br\'ezis and J.-L. Lions, eds.), vol.~2, Research Notes in Math., no.~60,
  Pitman, 1982, pp.~155--199.

\bibitem{GV85c}
\bysame, \emph{The global {C}auchy problem for the nonlinear {S}chr\"odinger
  equation revisited}, Ann. Inst. H. Poincar\'e Anal. Non Lin\'eaire \textbf{2}
  (1985), 309--327.

\bibitem{KPV00}
C.~Kenig, G.~Ponce, and L.~Vega, \emph{On the concentration of blow up
  solutions for the generalized {K}d{V} equation critical in {$L\sp 2$}},
  Nonlinear wave equations (Providence, RI, 1998), Contemp. Math., vol. 263,
  Amer. Math. Soc., Providence, RI, 2000, pp.~131--156.

\bibitem{KeraaniPhD}
S.~Keraani, \emph{{\'E}tudes de quelques r\'egimes asymptotiques de
  l'\'equation de {S}chr\"odinger}, Ph.D. thesis, Universit\'e Paris-Sud,
  Orsay, 2000.

\bibitem{Keraani01}
\bysame, \emph{On the defect of compactness for the {S}trichartz estimates of
  the {S}chr\"odinger equations}, J. Differential Equations \textbf{175}
  (2001), no.~2, 353--392. \MR{1 855 973}

\bibitem{Kwong}
Man~Kam Kwong, \emph{Uniqueness of positive solutions of ${\Delta} u-u+u\sp
  p=0$ in ${\R^n}$}, Arch. Rational Mech. Anal. \textbf{105} (1989), no.~3,
  243--266. \MR{90d:35015}

\bibitem{MerleDuke}
F.~Merle, \emph{Determination of blow-up solutions with minimal mass for
  nonlinear {S}chr\"odinger equations with critical power}, Duke Math. J.
  \textbf{69} (1993), no.~2, 427--454. \MR{94b:35262}

\bibitem{MerleRaphaelInv}
F.~Merle and P.~Rapha\"el, \emph{On universality of blow-up profile for ${L}^2$
  critical nonlinear {S}chr\"odinger equation}, Invent. Math. \textbf{156}
  (2004), 565--672.

\bibitem{MerleVega98}
F.~Merle and L.~Vega, \emph{Compactness at blow-up time for ${L}\sp 2$
  solutions of the critical nonlinear {S}chr\"odinger equation in 2{D}},
  Internat. Math. Res. Notices (1998), no.~8, 399--425. \MR{99d:35156}

\bibitem{MetivierSchochet98}
G.~M{\'e}tivier and S.~Schochet, \emph{Trilinear resonant interactions of
  semilinear hyperbolic waves}, Duke Math. J. \textbf{95} (1998), no.~2,
  241--304.

\bibitem{MVV}
A.~Moyua, A.~Vargas, and L.~Vega, \emph{Restriction theorems and maximal
  operators related to oscillatory integrals in {$\mathbb R\sp 3$}}, Duke Math.
  J. \textbf{96} (1999), no.~3, 547--574. \MR{2000b:42017}

\bibitem{Niederer}
U.~Niederer, \emph{The maximal kinematical invariance groups of {S}chr\"odinger
  equations with arbitrary potentials}, Helv. Phys. Acta \textbf{47} (1974),
  167--172. \MR{51 \#2511}

\bibitem{Perelman}
G.~Perelman, \emph{On the formation of singularities in solutions of the
  critical nonlinear {S}chr\"odinger equation}, Ann. Henri Poincar\'e
  \textbf{2} (2001), no.~4, 605--673.

\bibitem{Rauch91}
J.~Rauch, \emph{Partial differential equations}, Graduate Texts in Math., vol.
  128, Springer-Verlag, New York, 1991.

\bibitem{Strauss77}
W.~A. Strauss, \emph{Existence of solitary waves in higher dimensions}, Comm.
  Math. Phys. \textbf{55} (1977), no.~2, 149--162.

\bibitem{Strichartz}
R.~Strichartz, \emph{Restrictions of {F}ourier transforms to quadratic surfaces
  and decay of solutions of wave equations}, Duke Math. J. \textbf{44} (1977),
  no.~3, 705--714.

\bibitem{TaoUtah03}
T.~Tao, \emph{Recent progress on the {R}estriction conjecture}, {\tt
  arXiv:math.CA/0311181}, 2003, Lecture notes, Park City, Utah.

\bibitem{TaoGAFA03}
\bysame, \emph{A sharp bilinear restriction estimate on paraboloids}, Geom.
  Funct. Anal. \textbf{13} (2003), no.~6, 1359--1384.

\bibitem{Weinstein83}
M.~I. Weinstein, \emph{Nonlinear {S}chr\"odinger equations and sharp
  interpolation estimates}, Comm. Math. Phys. \textbf{87} (1982/83), no.~4,
  567--576. \MR{84d:35140}

\bibitem{Weinstein86}
\bysame, \emph{On the structure and formation of singularities in solutions to
  nonlinear dispersive evolution equations}, Comm. in Partial Diff. Eq.
  \textbf{11} (1986), no.~5, 545--565.

\bibitem{Yajima87}
K.~Yajima, \emph{Existence of solutions for {S}chr\"odinger evolution
  equations}, Comm. Math. Phys. \textbf{110} (1987), 415--426.

\end{thebibliography}

\end{document}